\documentclass{amsart}
\usepackage[nohug]{diagrams}
\usepackage[mathscr]{eucal}
\usepackage{amsfonts}
\usepackage{amsmath}
\usepackage{amsthm}
\usepackage{amssymb}
\usepackage{latexsym}

\begin{document}

\newcommand{\nc}{\newcommand}
\newcommand{\bom}{{_{\mathbf{\omega}}}}
\newcommand{\st}{\divideontimes}
\def\neweq{\setcounter{theorem}{0}}
\newtheorem{theorem}[]{Theorem}
\newtheorem{proposition}[]{Proposition}
\newtheorem{corollary}[]{Corollary}
\newtheorem{lemma}[]{Lemma}
\theoremstyle{definition}
\newtheorem{definition}[]{Definition}
\newtheorem{remark}[]{Remark}
\newtheorem{conjecture}[equation]{Conjecture}
\newcommand{\dis}{{\displaystyle}}
\def\question{\noindent\textbf{Question.} }
\def\remark{\noindent\textbf{Remark.} }
\def\proof{\medskip\noindent {\textsl{Proof.} \ }}
\def\endproof{\hfill$\square$\medskip}
\def\str{\rule[-.2cm]{0cm}{0.7cm}}
\newcommand{\beq}{\begin{equation}\label}
\newcommand{\aand}{\quad{\text{\textsl{and}}\quad}}
%
%{\underset{\rightarrow}{{\mathtt  lim}}}
%{\vphantom{x}\smash{\buildrel_\to\under x}}
%
\def\x{\stackrel{\rightarrow}{x}}
%{\vphantom{x}\smash{\buildrel^{_\rightarrow}\over x}}
\newcommand{\iso}{{\;\;\stackrel{\sim}{\longrightarrow}\;\;}}
\newcommand{\cd}{\!\cdot\!}
\newcommand{\vi}{${\sf {(i)}}\;$}
\newcommand{\vii}{${\sf {(ii)}}\;$}
\newcommand{\viii}{${\sf {(iii)}}\;$}
\newcommand{\viv}{${\sf {(iv)}}\;$}
\newcommand{\sset}{\subset}
\newcommand{\G}{\Gamma}
\newcommand{\tdr}{{\widetilde{\sf{DR}}}}
\newcommand{\id}{{{\mathtt {Id}}}}
\newcommand{\im}{{{\mathtt {Im}}}}
\newcommand{\Ind}{{{\mathtt {Ind}}}}
\newcommand{\sq}{{$\enspace\square$}}
\newcommand{\into}{\,\,\hookrightarrow\,\,}
\newcommand{\otni}{\,\,\hookleftarrow\,\,}
\newcommand{\too}{\,\,\longrightarrow\,\,}
\newcommand{\onto}{\,\,\twoheadrightarrow\,\,}
\newcommand{\bull}{^{^{_\bullet}}}
\newcommand{\sto}{\rightsquigarrow}
\newcommand{\ad}{{\mathtt{{ad}}}}
\newcommand{\Ad}{{\mathtt{{Ad}}}}
\newcommand{\Lie}{{\mathtt{Lie}^{\,}}}
\newcommand{\Spec}{{\mathtt{Spec}}}
\newcommand{\Hilb}{{\mathtt{Hilb}}}
\newcommand{\End}{{\mathtt{End}}}
\newcommand{\Endo}{{\mathtt{End}^{^{\mathsf{opp}}}}}
\newcommand{\Hom}{{\mathtt{Hom}}}
\newcommand{\Coh}{{\mathtt{Coh}}}
\newcommand{\Qcoh}{{\mathtt{Qcoh}}}
\newcommand{\Proj}{{\mathtt{Proj}}}
\newcommand{\Mat}{{\mathtt{Mat}}}
\newcommand{\Aut}{{\mathtt{Aut}}}
\newcommand{\DGAut}{{\mathtt{DGAut}}}
\newcommand{\Mod}{{\mathtt{Mod}}}
\newcommand{\MOD}{{\mathtt{MOD}}}
\newcommand{\Ideals}{{\mathtt{Ideals}}}
\newcommand{\DGMod}{{\mathtt{DGMod}}}
\newcommand{\Modi}{{\mathtt{Mod}}_{\infty}}
\newcommand{\Hi}{{\mathcal{H}}_{\infty}}
\newcommand{\Di}{{\mathcal{D}}_{\infty}}
\newcommand{\Com}{{\mathtt{Com}}}
\newcommand{\Ker}{{\mathtt{Ker}}}
\newcommand{\bX}{\,{\overline{\!X}}}
\newcommand{\tA}{{\tilde{A}}}
\newcommand{\tM}{{\tilde{M}}}
\newcommand{\tE}{{\tilde{E}}}
\newcommand{\tI}{{\tilde{I}}}
\newcommand{\tQ}{{\tilde{Q}}}
\newcommand{\tN}{{\tilde{N}}}
\newcommand{\bH}{{\underline{H\!}}\,}
\newcommand{\Qgr}{{\mathtt{Qgr}}}
\newcommand{\GrM}{{\mathtt{GrMod}}}
\newcommand{\grm}{{\mathtt{grmod}}}
\newcommand{\Coker}{{\mathtt{Coker}}}
\newcommand{\Span}{{\mathtt{span}}}
\newcommand{\Ann}{{\mathtt{Ann}}}
\renewcommand{\dim}{{\mathtt{dim}}}
\renewcommand{\deg}{{\mathtt{deg}}}
\newcommand{\Par}{{\mathtt{Part}}}
\newcommand{\co}{{\mathtt{cont}}}
\newcommand{\rk}{{\mathtt{rk}}}
\newcommand{\ch}{\chi}
\newcommand{\vr}{\boldsymbol{\varrho}}
\newcommand{\bHom}{\underline{\rm Hom}}
\newcommand{\bExt}{\underline{\rm Ext}}
\newcommand{\lrho}{\acute{\rho}}
\newcommand{\rrho}{\grave{\rho}}
\newcommand{\dr}{{{\mathsf {DR}}}}
\newcommand{\chom}{{\mathcal{H}om}}
\newcommand{\grd}{{\mathtt{gr}}}
\newcommand{\Ext}{{\rm{Ext}}}
\newcommand{\trace}{{\mathtt{Trace}}}
\newcommand{\Tr}{{{\mathtt{Tr}}}}
\newcommand{\tr}{{{\mathtt{Tr}}}}
\newcommand{\irr}{{{\mathsf {Irred}}}}
\newcommand{\rad}{^{^{\mathsf{rad}}}}
\newcommand{\hr}{{\mathfrak{h}^{^{\mathsf{reg}}}}}
\newcommand{\gr}{{\mathfrak{g}^{^{_{\mathsf{rs}}}}}}
\newcommand{\dd}{{\mathcal{D}}}
\newcommand{\GL}{\mathtt{GL}}
\newcommand{\SL}{\operatorname{SL}}
\newcommand{\SSp}{\mathtt{Sp}}
\newcommand{\hh}{{\mathsf{H}}}
\newcommand{\thh}{{\mathsf{\tilde{H}}}}
\newcommand{\ehe}{{\mathbf{e}\!\mathsf{H}\!\mathbf{e}}}
\newcommand{\eve}{{\mathbf{e}\mathcal{H}\mathbf{e}}}
\newcommand{\e}{{\mathbf{e}}}
\newcommand{\cc}{{\mathbf{c}}}
\newcommand{\ka}{\kappa}
\newcommand{\MM}{{\mathcal{M}}}
\newcommand{\RR}{{\mathcal{R}}}
\newcommand{\AB}{{\boldsymbol{A}}}
\newcommand{\BA}{{\boldsymbol{B}}}
\newcommand{\LB}{{\boldsymbol{L}}}
\newcommand{\sym}{{\mathsf{Sym}}}
\newcommand{\dertr}{{\mathtt{Der}}_{\!_{\mathsf {tr}}}}
\newcommand{\difp}{\mathtt{Diff}_{_{\!{\mathscr{P}}}}}
\newcommand{\derp}{\mathtt{Der}_{_{\!{\mathscr{P}}}}}
\newcommand{\sgn}{{\mathsf{sign}}}
\newcommand{\Ug}{{{\cal{U}}\g}}
\newcommand{\Uh}{{{\cal{U}}\h}}
\newcommand{\triv}{{\mathtt {triv}}}
\newcommand{\sign}{{\mathtt {sign}}}
\newcommand{\sch}{\boldsymbol{\mathsf{s}}}
\newcommand{\mon}{\boldsymbol{\mathsf{m}}}
\newcommand{\aut}{{\mathtt {Aut}}}
\newcommand{\tdash}{\mbox{\tiny{-}}}
\newcommand{\Th}{\Theta}
\newcommand{\lL}{\underround}
\newcommand{\rR}{{}}
\newcommand{\bphi}{\bar{\phi}}
\newcommand{\ee}{{\mathfrak{e}}}
\newcommand{\h}{{\mathfrak{h}}}
\newcommand{\cm}{{\mathbb{M}}}
\newcommand{\var}{\kappa}
\newcommand{\dvar}{{\mathsf{D}}_\kappa}
\newcommand{\Eu}{{{\mathtt {Eu}}}}
\newcommand{\om}{\Omega}
\newcommand{\A}{{\mathsf{A}}_{\infty}}
\newcommand{\Om}{\Omega}
\newcommand{\VV}{\mathsf{Vect}}
\newcommand{\wh}{\widehat}
\newcommand{\bh}{\bar{H}}
\newcommand{\la}{\label}
\newcommand{\greg}{{\mathfrak{g}^{^{_{\mathsf{reg}}}}}}
\newcommand{\hreg}{{\mathfrak{h}^{^{_{\mathsf{reg}}}}}}
\newcommand{\V}{{\mathscr{V}}_\var}
\newcommand{\eij}{{\epsilon_{_{ij}}}}
\newcommand{\eps}{{\epsilon}}
\newcommand{\bul}{^\bullet\!}
\newcommand{\bulp}{^\bullet_{_{\!\PP\!}}}
\newcommand{\hv}{{\mathcal{H}}}
\newcommand{\va}{\varkappa}
\def\bs#1{\boldsymbol{#1}}
\def\ms#1{\mathcal{#1}}
\def\C{{\mathbb{C}}}
\def\X{{\bar{X}}}
\def\Y{{\bar{Y}}}
\def\j{{\bar{j}}}
\def\g{{\bar{g}}}
\def\Q{{\mathbb{Q}}}
\def\R{{\mathsf{R}}}
\def\qq{{\mathcal{Q}}}
\def\ms#1{\mathcal{#1}}
\def\AA{{\mathbb{A}}}
\def\II{{\mathbb{I}}}
\def\rr{{\mathcal{R}}}
\def\Z{{\mathbb{Z}}}
\def\oo{{\mathcal O}}
\def\K{\boldsymbol{K}}
\def\tx{\tilde{x}}
\def\ty{\tilde{y}}
\def\ga{\mathfrak{a}}
\def\omeg{\mathfrak{w}}
\def\bi{\boldsymbol{i}}
\def\KK{ \tilde{\boldsymbol{K}}}
\def\L{\boldsymbol{L}}
\def\M{\boldsymbol{M}}
\def\N{\boldsymbol{N}}
\def\xx{\boldsymbol{x}}
\def\E{\boldsymbol{E}}
\def\U{{\mathcal U}}
\def\ccr{\C[R]^W_{_{\sf reg}}}
\def\ccz{\Z[R]^W}
\def\up{{\mathcal U}^{\PP\!}\!}
\def\upa{{{\mathcal U}^{\PP\!}\!A}}
\def\pp{_{_{\!\PP}}}
\def\ll{{\mathcal L}}
\def\LL{{\mathbf{L}}}
\def\D{{\mathcal{D}}}
\def\rep{{\mathsf{Rep}}}
\def\Pf{{\it Proof}}
\def\PP{\mathbb{P}}
\def\O{{\sf O}}
\def\aa{{\mathcal{A}}}
\def\eu{{\mathsf{eu}}}
\def\reg{{\!}^{^{\mathsf{reg}}}}
\def\bb{{\mathcal{B}}}
\def\T{{\mathsf T}}
\def\sll2{{\mathfrak{s}\mathfrak{l}}_2}
\def\gln{{\mathfrak{g}\mathfrak{l}}_n(\C)}
\def\gl{{\mathfrak{g}\mathfrak{l}}}
\def\sln{{\mathfrak{s}\mathfrak{l}}_n(\C)}
\def\glinfty{{\mathfrak{g}\mathfrak{l}}_\infty}
\def\Cr{\C^{\mathsf{reg}}}
\def\CC{{\mathcal{C}}}
\def\HH{{\mathcal{H}}}
\def\F{{\mathcal{F}}}
\def\FF{{\mathsf{F}}}
\def\s{{\mathbb{S}}}
\def\sn{{{\mathbb{S}}_n}}
\def\w{{\mathsf{W}}}
\def\bpa{{{\mathsf{B}}_{\bullet}^\PP\!A}}
\def\ccirc{{{}_{^{^\circ}}}}
\def\dertr{{\der_{\!_{\om\otimes\tr}}}}
\newcommand{\hdv}{{{\vphantom{H}
\smash{\buildrel_{\,_{_{_{\bullet\bullet}}}}\over {\sf H}}}}_{\!\kappa}}
\numberwithin{equation}{section}
\title{$\A$-Modules and Calogero-Moser Spaces}
\author{Yuri Berest}
\address{Department of Mathematics,
Cornell University, Ithaca, NY 14853-4201, USA}
\email{berest@math.cornell.edu}
\thanks{Berest's work partially supported by NSF grant DMS 04-07502 and
A.~P.~Sloan Research Fellowship.}
\author{Oleg Chalykh}
\address{Department of Mathematics, University of Leeds, Leeds LS2 9JT, UK}
\email{oleg@maths.leeds.ac.uk}
\maketitle

\section{Introduction}
The Hilbert schemes $ \Hilb_n(\C^2) $ of points on $ \C^2 $ have a rich geometric structure
with many interesting links to representation theory, combinatorics and integrable systems.
One reason for this is perhaps that the points of $ \Hilb_n(\C^2) $ admit a few different
algebraic incarnations which underlie the geometric properties of $ \Hilb_n(\C^2)$.
Specifically, the space $\,\Hilb(\C^2) := \bigsqcup_{n \geq 0} \Hilb_n(\C^2) \,$ parametrizes
\begin{enumerate}
\item the ideals of finite codimension in the polynomial algebra $ A_0 := \C[x,y] \,$;
\item the isomorphism classes of finite-dimensional representations $\, (V, \,\bar{i}) \,$
      of $ A_0 $ with a fixed cyclic vector $ \bar{i} \in V $;
\item the isomorphism classes of finitely generated rank $1$ torsion-free $A_0$-modules;
\item the isomorphism classes of rank $1$ torsion-free coherent sheaves on $ \PP^2(\C) $
      ``framed'' over the line at infinity.
\end{enumerate}
The relations between these objects are well known and almost immediate. Thus, (1) is essentially
a definition of (closed) points of  $\,\Hilb(\C^2) $. The bijection (1) $ \to $ (2) is given by
taking the quotient $\, M \mapsto A_0/M \,$ modulo a given ideal and letting $\, \bar{i} \,$ be
the image of $ 1 \in A_0 $ in $ A_0/M $. The inverse map (2) $\to $ (1) is then
defined by assigning to a given cyclic module its annihilator in $ A_0$.
The correspondence (1) $\leftrightarrow$ (3) follows from the
fact that every f.~g. rank $1$ torsion-free $A_0$-module is isomorphic
to a unique ideal of finite codimension in $ A_0 $. Finally, the bijection (3) $ \to $ (4)
can be constructed geometrically by extending $ A_0$-modules to coherent sheaves on $ \PP^2 $,
and its inverse by restricting such sheaves via the natural embedding $\, \C^2 \into \PP^2 \,$.

Now, let us ``quantize'' the affine plane $\, \C^2 \,$ replacing the commutative polynomial ring $ A_0 $
by the first complex Weyl algebra $\, A_1 := \C \langle x,y \rangle /(xy-yx-1)\,$. One can ask
then the natural (though, perhaps, very na\"{\i}ve) question:
{\it What happens to the above bijections?}
At first glance, this question does not make sense since only (3) has a clear analogue for the Weyl algebra.
However, following an idea of Le Bruyn \cite{LeB}, we can replace $ \PP^2 $ (or rather, the category
$ \Coh(\PP^2) $ of coherent sheaves on $ \PP^2 $) by a quantum projective plane
$\, \PP^2_{\!q}\,$ and identify a class of objects in $ \Coh(\PP^2_{\!q}) $
that are natural analogues (deformations) of (4).
As a result, we can extend the bijection (3) $\leftrightarrow$ (4) to the noncommutative
case (see \cite{BW2}).

In this paper we make one step further suggesting what might be a ``quantum analogue'' of a
finite-dimensional cyclic representation of $ A_0 $. Our main observation is that
the Weyl algebra $ A_1 $ does have finite-dimensional modules $ V $,
which can be related to its ideals in an essentially canonical way, provided we
relax the associativity assumption on the action of $ A_1 $, i.~e. assume that
$$
(v.a).b \not= v.(ab) \quad \mbox{for some}\ \, a,b \in A_1 \ \mbox{and}\ v \in V \ .
$$
As we will see, such ``non-associative representations'' of $ A_1 $ have a natural origin
from the point of view of deformation theory. To define them we should think of $ A_1 $ not
as an associative algebra but as an {\it $\A$-algebra}, and thus work not with
(complexes of) $ A_1$-modules but with {\it $\A$-modules} over $A_1$.

\def\i{{\bar{i}}}

To explain this idea we return for a moment to the commutative case. By definition,
a cyclic representation of $ A_0 $ is an $A_0$-module $ V $ generated by a single
vector $ \i \in V\,$. Giving a pair $ (V, \, \i) $ is then equivalent to giving a surjective
$A_0$-linear map $\, A_0 \to V $, $\,1 \mapsto \i \,$, which, in turn, can be written as
a two-term complex of $A_0$-modules
\begin{equation}
\la{2com}
 0 \to A_0 \to V \to 0 \ .
\end{equation}

Now, for any associative algebra $ A $ there is a natural (``interpretation'')
functor $ \Com(A) \to \Modi(A) $ from the category of complexes of $A$-modules
to that of $\A$-modules\footnote{We will review the definition and basic properties of
$\A$-modules in Section~\ref{Sec2}.} over $A$. This functor is faithful, but neither full nor
surjective: in other words, $ \Com(A) $ can be viewed as a subcategory of $ \Modi(A) \,$,
but $ \Modi(A) $ has more objects and more morphisms than $ \Com(A)$.

If we deform now $ A_0 $ to $ A_1 $ via the family of algebras $\,
A_{\hbar} := \C \langle x,y \rangle / (xy - yx - \hbar) $, the complex \eqref{2com} with
$\,  0 < \dim(V) < \infty \,$  does not admit deformations in $ \Com(A_{\hbar}) \,$
(as $ A_{\hbar} $ has no non-trivial finite-dimensional modules except for $ \hbar = 0 $).
However, it can be deformed naturally within the larger categories $ \Modi(A_{\hbar}) $.
The resulting $\A$-module $ \K  $ can still be represented by a two-term complex of vector spaces
$\, 0 \to K^0 \to K^1 \to 0 \,$, with $ K^1 $ being finite-dimensional, but the action of
$ A_\hbar $ on $ \K $ will not be strictly associative. Letting $ \hbar = 1 $ and
restricting to $ K^1 $, we get thus a finite-dimensional ``non-associative representation'' of $A_1$.
We will characterize such representations (or rather, the corresponding $\A$-modules $\K $)
axiomatically and relate them to the rank $1$ torsion-free right modules (ideals) of $ A_1 $.

In the commutative case, the ideal (class) of $ A_0 $ corresponding to a cyclic representation
$ (V, \, \i) $ is determined by cohomology of the complex \eqref{2com}. For the Weyl algebra,
the relation is now similar: every ideal $ M $ of $A_1$ embeds in the corresponding
$ \K $ as $ \A$-module, and this embedding is a quasi-isomorphism in $ \Modi(A_1)$.
Thus, relative to $M$, the $\A$-module $ \K $ plays the role of a certain resolution
in $ \Modi(A_1)$ whose properties resemble the properties of minimal resolutions (envelopes)
in classical homological algebra. We will therefore refer to $\K$ as an {\it $\A$-envelope} of $M$.

In view of non-associativity, the action of $ x $ and $ y $ of $ A_1 $ on the $\A$-module
$\, \K = K^0 \oplus K^1 \,$ is not subject to the canonical commutation relation.
Instead, when restricted to $ K^1 $, the corresponding endomorphisms $\X $ and $ \Y $ satisfy the
``rank-one'' condition: $\, \rk([\X, \,\Y] + \id) = 1 \,$. We will show that $ \K $ can be
uniquely reconstructed from the data $ (K^1, \X, \Y) $ up to strict isomorphism.
Thus we establish a bijection  between the set $ \MM $ of strict isomorphism classes of $\A$-envelopes and the disjoint union $ \CC  $ of the Calogero-Moser varieties $\CC_n \,$ (see the definition below).
On the other hand, an object of $ \Modi(A_1) $ satisfying the axioms of $\A$-envelopes is uniquely determined by its cohomology which, in turn, is given by a rank $1$ torsion-free $A_1$-module.
Hence, we have also a bijection $ \MM \leftrightarrow \RR \,$, where
$ \RR $ is the set of isomorphism classes of (right) ideals of $A_1$. Combining these last two
bijections, we arrive at the {\it Calogero-Moser correspondence}
$\, \RR \leftrightarrow \CC \,$, which gives a geometric classification of ideals of $ A_1 $.

The correspondence $ \RR \leftrightarrow \CC $ was first proved in \cite{BW1} by combining some earlier results of Cannings-Holland \cite{CH} and Wilson~\cite{W}.
Two other proofs using the methods of noncommutative projective geometry and representation theory of quivers can be found in \cite{BW2} and in the appendix to \cite{BW2}.
All three proofs are fairly involved and indirect, especially in contrast with elementary
arguments in the commutative case. A proof given in this paper results from our attempt to extend those
arguments to the noncommutative case. As an indication of this attempt being worth-while, we mention a simple formula for the Calogero-Moser map $\, \omega: \CC \to \RR \,$, which appears naturally in our
approach but seems to be missing (or implicit) in earlier papers\footnote{Actually, this formula is a
``noncommutative version'' of a remarkable formula of G.~Wilson for the rational Baker function of the
KP hierarchy (see \cite{W}). It can be deduced by comparing the results of \cite{BW2} and \cite{W}
(see Notes in \cite{BW3}, p.~116).}.

First, we recall that the variety $ \CC_n $ can be defined as a quotient of
the space of matrices $\, \{(\X, \Y, \i, \j) \,:\, \X, \Y \in \End(\C^n) \, , \, \i \in
\Hom(\C,\C^n) , \, \j \in \Hom(\C^n, \C)\} \,$ satisfying the equation
$\, [\X,\,\Y] + \id_n = \i\,\j \,$ modulo a natural action of $ \GL_n(\C) \,$ (see \cite{W}).
Now, given a point  $\, [(\X, \Y, \i, \j)] \,$ of $ \CC_n $, we claim that the class of
$ \RR $ corresponding to it under the bijection $ \omega $ can be represented by the
(fractional) ideal
$$
M =  \kappa \,\det(\X - x)\,A_1 + \det(\Y - y)\,A_1\ ,
$$
where $ \kappa  $ is given by the expression $\, 1 - \j\,(\Y - y)^{-1} (\X-x)^{-1}\i \,$
in the quotient skew-field of $A_1$. Surprisingly, in the commutative case, there seems to be
no analogue of such an explicit presentation of ideals.

A few words about the organization of this paper: it consists of nine sections, each starting with
a brief introduction. There is also an appendix containing an alternative (geometric) construction
of $\A$-envelopes. In the last section we discuss the question of functoriality of the Calogero-Moser
correspondence which was originally our motivation for the present work. As often happens, we have not
clarified it completely, but we hope some more details will appear elsewhere.

\subsection*{Acknowledgements}
{\footnotesize
We thank V.~Ginzburg, D.~Kazhdan, B.~Keller, M.~Kontsevich, A.~Rosenberg, P.~Smith, T.~Stafford,
B.~Tsygan, M.~Van den Bergh and G.~Wilson for many interesting discussions, questions and comments.
We are especially indebted to B.~Keller for advice in the early stages of the present work,
and to G.~Wilson for his encouragement and gentle criticism.
The first author is grateful to the Mittag-Leffler Institut for its hospitality
during the period May-June 2004, where the main part of this paper was written. }

\section{$\A$-Modules and Morphisms}
\la{Sec2}
In this section we review the definition of $\A$-modules and their homomorphisms.
These concepts can be defined naturally over an arbitrary $\A$-algebra (see \cite{Ka},
\cite{K1}). However, in the present paper we deal mostly with usual associative algebras
and thus we restrict our discussion below to this special case.

\subsection{$\A$-modules}

Let $ A $ be a unital associative algebra over a field $ k\,$. In what follows
we will often think of $A$ as being $ \Z$-graded with single nonzero component
$\, A^0 = A \,$ in degree $ 0\,$.

A (right) {\it $ \A$-module} over $ A $ is a $ \Z$-graded $k$-vector space $ \K =
\bigoplus_{p \in \Z} K^p $ equipped with a sequence of homogeneous multilinear operations
$$
m_n:\, \K \otimes A^{\otimes (n-1)} \to \K \ , \quad n\geq 1 \ .
$$
These operations are subject to the following conditions.

First, $\, m_1: \K \to \K \,$ has degree $ +1 $ and satisfies the equation
\begin{equation}
\label{1.1}
(m_1)^2 = 0 \ .
\end{equation}
Thus, $\, \K \, $ is a complex of vector spaces with differential $ m_1\,$.

Second, $\, m_2: \K \otimes A \to \K \,$ has degree $ 0 $ and commutes with $ m_1 \,$:
\begin{equation}
\label{1.2}
m_1(m_2(x,a)) = m_2(m_1(x), a) \ , \quad x\in \K\ , \ a\in A \ .
\end{equation}
Thus, $ m_2 $ may be thought of as an action of $ A $ on the complex
$\, (\K, m_1) \,$. This action, however, need not be associative. The
corresponding associativity diagram
\[
\begin{diagram}
\K \otimes A \otimes A  & \rTo^{m_2 \otimes \id \,} & \K \otimes A \\
\dTo^{\id \otimes m_A} & \rdTo^{m_3} & \dTo_{m_2} \\
\K \otimes A & \rTo^{m_2} &  \K \\
\end{diagram}
\]
commutes only ``up to homotopy,'' which is specified by the next operation $ m_3\,$.

Thus, $\, m_3: \K \otimes A \otimes A \to \K \,$ is a map of degree $ -1 $
satisfying
\begin{equation}
\label{1.3}
m_2(m_2(x,a),b) - m_2(x, ab) = m_3(m_1(x),a,b) + m_1(m_3(x,a,b))\
\end{equation}
for all $ x \in \K $ and $ a,b \in A\,$.

In general, the maps $ m_n $ have degree $ 2-n $ and satisfy the following
algebraic relations (called the {\it strong homotopy relations}):
\begin{eqnarray}
\label{1.4}
\lefteqn{\sum_{i=1}^{n-2} (-1)^{n-i} m_{n-1}(x, a_1, a_2, \ldots, a_i a_{i+1}, \ldots, a_{n-1})
+} \\
&&\sum_{j=1}^{n}(-1)^{n-j} m_{n-j+1}(m_j(x, a_1, a_2, \ldots, a_{j-1}), a_j,
\ldots, a_{n-1}) = 0\ . \nonumber
\end{eqnarray}
for all $\, x \in \K \,$ and $\, a_1, a_2, \ldots, a_{n-1} \in A\,$.

Since $ A $ is a unital algebra, it is natural to work with {\it unital}\, $\A$-modules: thus,
in addition to (\ref{1.4}), we will assume that $\, m_2(x,1) = x \,$
and $\, m_n(x,\ldots, 1, \ldots) = 0 \,$, $\, n \geq 3\,$, for all $ x \in \K \,$.

Observe that if $ m_n \equiv 0 $ for all $ n \geq 3 $ then $ m_2 $ is associative,
and $ (\K, m_1, m_2) $ can be identified with a usual complex of (right) $A$-modules.
Moreover, if $\K$ has only finitely many (say, $N$) nonzero components, all being
in non-negative degrees, then we have
$\, m_n \equiv 0 \,$ for $\, n > N+1 \,$, because $ \,\deg(m_n) = 2-n\,$.
This does not mean, however, that any choice of
linear maps $ (m_1, m_2, \ldots, m_{N}) $ satisfying the first $N$ equations of
(\ref{1.4}) extends to an $\A$-structure on $ \K \,$. In general, the higher homotopy
relations impose certain obstructions; we will need the following easy result showing
that no such obstructions arise in the special case of complexes with two components.
\begin{lemma}
\la{L1}
Let $ \K := [\,0 \to K^0 \stackrel{m_1}{\longrightarrow} K^1 \to 0 \,] $
be a two-term complex of vector spaces equipped with a surjective differential
$ m_1 $ and operations $ m_2 $ and $ m_3 $ satisfying  \eqref{1.2} and \eqref{1.3}.
Then the triple $ (m_1, m_2, m_3) $ extends to a (unique) structure of $\A$-module
on $ \K \,$.
\end{lemma}
\begin{proof}
Since $ \K $ has nonzero components only in degrees $ 0 $ and $ 1 $
we have $\, m_n \equiv 0 $ for all $ n \geq 4\,$. It remains to check that
the sequence of maps $ (m_1, m_2, m_3, 0, 0, \ldots) $ satisfies the relations (\ref{1.4}).
These relations hold automatically for $ n > 4 $, while for $ n = 4 $ we have
the apparent compatibility condition:
\begin{equation}
\la{1.10}
-m_3(x, ab, c) + m_3(x,a,bc) + m_3(m_2(x,a),b,c) - m_2(m_3(x,a,b),c) = 0\ .
\end{equation}
Letting $\, x = (u,v) \in K^0 \oplus K^1 \,$, we may rewrite
(\ref{1.2}), (\ref{1.3}) and (\ref{1.10}) as
\begin{equation}
\la{1.2e}
m_1(m_2^0(u,a)) = m_2^1(m_1(u), a)\ ,
\end{equation}
\begin{equation}
\la{1.3e}
m_2^0(m_2^0(u,a),b) - m_2^0(u, ab) = m_3(m_1(u), a,b) \ ,
\end{equation}
\begin{equation}
\la{1.31e}
m_2^1(m_2^1(v,a),b) - m_2^1(v, ab) = m_1(m_3(v, a,b))\ ,
\end{equation}
\begin{equation}
\la{1.10e}
-m_3(v, ab, c) + m_3(v,a,bc) + m_3(m_2^1(v,a),b,c) - m_2^0(m_3(v,a,b),c) = 0\ ,
\end{equation}
where $\, m_2^0: K^0 \otimes A \to K^0 \,$ and $\, m_2^1: K^1 \otimes A \to K^1 \,$
denote the two nontrivial components of $ m_2 \,$.
Since $\, m_1: K^0 \to K^1 \,$ is surjective,  (\ref{1.3e}) uniquely determines
$ m_3 $ in terms of $ m_1 $ and $ m_2 \,$, and (\ref{1.10e}) is easily seen to be an
algebraic consequence of (\ref{1.2e}) and (\ref{1.3e}).
\end{proof}

\subsection{Morphisms of $\A$-modules}
\la{S2.2}
A {\it morphism} $\, f: \K \to \L \,$ between two $\A$-modules over $ A $
is defined by a sequence of homogeneous linear maps
$$
f_n:\, \K \otimes A^{\otimes (n-1)} \to \L \ , \quad n\geq 1 \ ,
$$
which are subject to the following conditions.

First, $\, f_1: \K \to \L \,$ has degree $ 0 $ and commutes with differentials
on $ \K $ and $ \L \,$:
\begin{equation}
\label{1.5}
m_{1}^{L} \, f_1 = f_1 \, m_{1}^{K} \ .
\end{equation}
Thus, $ f_1 $ is a morphism of complexes of vector spaces $ (\K, m_{1}^{K}) $ and
$ (\L, m_{1}^{L})$. This, however, need not be $ A$-linear with respect to $ m_{2}^{L} $
and $ m_{2}^{K} \,$. The corresponding linearity diagram
\[
\begin{diagram}
\K \otimes A  & \rTo^{m_{2}^{K}} & \K \\
\dTo^{f_1 \otimes \id} & \rdTo^{f_2} & \dTo_{f_1} \\
\L \otimes A & \rTo^{m_{2}^{L}} &  \L \\
\end{diagram}
\]
commutes only ``up to homotopy'' to be specified by the next component of $ f\,$.

Thus, $\, f_2: \K \otimes A \to \L \,$ is a map of degree $ -1 $ satisfying
\begin{equation}
\label{1.6}
f_1(m_{2}^{K}(x,a)) - m_{2}^{L}(f_1(x), a) = f_2(m_{1}^{K}(x),a) + m_{1}^{L}(f_2(x,a))\
\end{equation}
for all $ x \in \K $ and $ a \in A\,$.

In general, the maps $ f_n $ have degree $\, 1 - n \,$ and satisfy
some infinite system of algebraic relations similar to \eqref{1.4}
(see \cite{K2}, (6.9)).

If both $ \K $ and $ \L $ have at most $ N $ nonzero components
(located in non-negative degrees), then $ f_n \equiv 0 $ for all $
n > N+1 \,$. Not any pair of linear maps $ (f_1, f_2) $ satisfying
\eqref{1.5} and \eqref{1.6} extends, in general, to an $
\A$-morphism. However, as in the case of structure maps (cf.
Lemma~1), the following result shows that no obstructions arise
for extending morphisms between two-term complexes.
\begin{lemma}
\la{L2}
Let $ \,\K \,$ and $ \,\L \,$ be $\,\A$-modules
having nonzero components only in degrees $ 0 $ and $ 1 \,$.
Assume that $\,  m_{1}^{K}: K^0 \to K^1 \,$ is surjective.
Then any pair of linear maps $(f_1, f_2)$ satisfying
\eqref{1.5} and \eqref{1.6} extends to a unique morphism $ f: \K \to \L \,$
of $\A$-modules.
\end{lemma}
\begin{proof}
The uniqueness is obvious, since we have $\, f_n \equiv 0 $ for
$ \, n \geq 3 \,$ by degree considerations. We need only to check that the sequence of maps
$ (f_1, f_2, 0, 0, \ldots) $  satisfies the higher homotopy relations, provided its first
two components satisfy (\ref{1.5}) and (\ref{1.6}). For $ n \geq 4 $, these relations hold trivially,
while for $ n = 3 $ we get the compatibility condition (cf. \cite{K2}, (6.9), $n=3$):
\begin{equation}
\la{1.11}
m_2(f_2(x,a),b) + m_3(f_1(x),a,b) = f_{2}(x,ab) - f_2(m_2(x,a),b) + f_1(m_3(x,a,b))\ . 
\end{equation}
As in the proof of Lemma~\ref{L1}, letting $ x = (u,v) \in K^0 \oplus K^1 $
we rewrite (\ref{1.5}), (\ref{1.6}) and (\ref{1.11}) in the form
\begin{equation}
\label{1.5e}
m_{1}(f_1^0(u)) = f_1^1(m_{1}(u)) \ ,
\end{equation}
\begin{equation}
\label{1.6e}
f_1^0(m_2^0(u,a)) - m_2^0(f_1^0(u),a) = f_2 (m_{1}(u),a) \ ,
\end{equation}
\begin{equation}
\label{1.61e}
f_1^1(m_2^1(v,a)) - m_2^1(f_1^1(v),a) = m_1(f_2(v,a)) \ ,
\end{equation}
\begin{equation}
\la{1.11e}
m_{2}^0(f_2(v,a),b) + m_{3}(f_1^1(v),a,b) = f_{2}(v,ab) - f_2(m_{2}^1(v,a),b) + f_{1}^0(m_3(v,a,b))\ . 
\end{equation}
Here $\, f_{1}^0: K^0 \to L^0 \,$ and $\, f_{1}^1: K^1 \to L^1 \,$ denote
the two components of the map $ f_1 \,$, and  $ f_2 $ is identified with
its only nonzero component $\, f_2: K^1 \otimes A \to L^0 \,$.
Since $ m_1: K^0 \to K^1 $ is surjective, equation (\ref{1.6e}) determines
$ f_2 $ in terms of $ f_1 \,$, and (\ref{1.61e}) is then an immediate
consequence of (\ref{1.5e}) and (\ref{1.6e}). Furthermore, applying
$\, m_2^0(\,\mbox{--}\,,b)\,$ to both sides of (\ref{1.6e}) and using (\ref{1.3e}) and
(\ref{1.5e}) we get after some trivial algebraic manipulations
\begin{eqnarray}
\la{1.111e}
\lefteqn{m_{2}^0(f_2(m_1(u),a),b) + m_{3}(f_1^1(m_1(u)),a,b) = } \\*[1ex]
&& f_{2}(m_1(u),ab) - f_2(m_{2}^1(m_1(u),a),b) + f_{1}^0(m_3(m_1(u),a,b))\ . \nonumber
\end{eqnarray}
Again, in view of surjectivity of $ m_1 \,$, (\ref{1.111e}) is equivalent to (\ref{1.11e}).
\end{proof}

The $\A$-morphisms $\, f: \K \to \L \,$ with $\,f_n \equiv 0\,$ for all $\, n \geq 2 \,$
are called {\it strict}. In view of \eqref{1.6}, $\, f \,$ being strict implies
that its first component $\, f_1 \,$ is $A$-linear. Thus, if $ \K $ and $ \L $ are usual complexes
of $A$-modules, strict $\A$-morphisms $\,\K \to \L \,$ can be identified with
usual morphisms of complexes. In general, working with arbitrary $\A$-modules, we will always
assume the identity morphisms to be strict.

\subsection{The category of $\A$-modules}

The (right unital) $ \A$-modules over $ A $ with (nonstrict) $\A$-morphisms form a category
which we denote $ \Modi(A)\,$. Since the usual complexes of modules over $A$
can be regarded as $\A$-modules (with higher operations $ m_n$, $\,n\geq3\,$, vanishing)
and the usual maps of such complexes can be identified with strict $\A$-morphisms,
the category $ \Com(A) $ can be interpreted as a subcategory of $ \Modi(A) \,$.
Note, however, though faithful, such an ``interpretation'' functor
$\, \Upsilon: \Com(A) \to \Modi(A) \,$ is neither full nor surjective: the category
$ \Modi(A) $ has more objects and more morphisms than $ \Com(A)\,$.

Assigning to an $\A$-module $ \K $ its {\it cohomology} $\, H^{n}(\K) \,$
(with respect to the differential $ m_1 $) and to an $\A$-morphism $\, f: \K \to \L \,$
the map $\, H^n(f) := H^n(f_1):\, H^n(\K) \to H^n(\L) \,$ induced on cohomology by
its first component gives a functor $\, H^n: \Modi(A) \to \Mod(A) \,$ with values
in the category of $A$-modules.
This functor is well defined, since each space $ H^n(\K) $ comes equipped with an action
of $ A $ induced by $ m_2 $, which is associative due to the homotopy relation (\ref{1.3}),
and each map $ f_1 $ is $A$-linear at the level of cohomology due to (\ref{1.6}).

We call a morphism $\, f: \K \to \L \,$  a {\it quasi-isomorphism} in $ \Modi(A) $
(in short, an {\it $\A$-quasi-isomorphism}) if the maps $\, H^n(f): H^{n}(\K) \to H^{n}(\L)\,$
are isomorphisms in $ \Mod(A) $ for all $ n \in \Z $. As in the classical case, the {\it derived category}
$\, \Di(A) \,$ of $\A$-modules can now be defined by universally localizing $\, \Mod(A) \,$ at the class of
all $\A$-quasi-isomorphisms. This notion, however, turns out to be ``redundant'' as the following important
result, due to Keller (see \cite{K1}), shows.
\begin{theorem}
\la{T1}
The canonical functor $\, \Upsilon: \Com(A) \to \Modi(A) \,$ descends to an embedding
$\, \D(\Upsilon): \D(A) \to \Di(A)\,$, which is an equivalence of (triangulated)
categories.
\end{theorem}
\begin{remark}
In \cite{K1} the category $ \Di(A) $ includes nonunital modules and thus,
strictly speaking, it is larger than the one we introduced above. In this nonunital setting
the functor $ \D(\Upsilon) $ is fully faithful but not surjective: the (essential) image
of $ \D(\Upsilon) $  consists of $\A$-modules which are unital at the cohomology level.
\end{remark}

\section{$\A$-Envelopes}
\la{Env}
Theorem~\ref{T1} shows that passing from usual (complexes of) modules
to $\A$-modules over $ A $ does not yield new quasi-isomorphism classes.
However, since $ \Di(A) $ has more objects than $\D(A) \,$, this does yield
{\it new representatives} of such classes. Being $\A$-modules, such representatives
come equipped with higher homotopy products, and these can be used to construct
new algebraic invariants of $A$-modules.

In this section we illustrate this general principle by looking at (probably) the
simplest nontrivial example: the rank one torsion-free modules over the Weyl algebra
$ A_1 \,$. Such modules are isomorphic to ideals of $ A_1\,$ and hence are all
projective (but not free). The classical (abelian) homological algebra
fails to produce any invariants that would allow one to distinguish such modules up
to isomorphism. However, as we will see below, such invariants --- the ``points'' of
the  Calogero-Moser varieties --- can be introduced via certain $\A$-modules
representing ideals in $ \Di(A_1) $. The properties of these $\A$-modules somewhat
resemble the properties of minimal resolutions (injective envelopes), and thus we
term them the $\A$-envelopes of our ideals.

\subsection{Axioms}
From now on, we assume $ k $ to be an algebraically closed field of characteristic zero, and
let $ A = A_1(k) $ denote the first Weyl algebra over $k\,$. We fix, once and for all,
two canonical generators $ x $ and $ y $ of $ A $ satisfying $\, xy - yx  = 1\,$,
and thus we distinguish two polynomial subalgebras  $ k[x] $ and $ k[y] $ in $A\,$.

Let $M$ be a rank $1$ finitely generated torsion-free module over $A\,$.
Using the canonical embedding $\, \Mod(A) \to \Modi(A) $, we will regard $ M $ as an object
of $ \, \Modi(A) $ (so that $ m_n^M \equiv 0 $ for $ n \not= 2 $ and $ m_2^M $ is the given
action of $A$ on $M$).
\begin{definition}
\la{D1}
An {\it $\A$-envelope}\, of $ M $ is an $\A$-quasi-isomorphism $\, r: M \to \K \,$,
where  $\, \K = K^0 \oplus K^1 \, $ is a unital $\A$-module over $ A $
with two nonzero components (in degrees $ 0 $ and $1$) and the structure maps
\begin{eqnarray}
&& m_1: \K \to \K\ , \quad (u,v) \mapsto (0, m_1(u))\ ,\nonumber\\
&& m_2: \K \otimes A \to \K\ ,\quad (u,v) \otimes a \mapsto (m_2^0(u,a), m_2^1(v,a))\ ,\nonumber\\
&& m_3: \K\otimes A \otimes A \to \K\ , \quad (u,v)\otimes a \otimes b
\mapsto (m_3(v,a,b),0)\ ,\nonumber
\end{eqnarray}
satisfying the axioms:
\begin{itemize}
\item {\it Finiteness:}
\begin{equation}
\la{2.1}
\dim_{k}\,K^1  < \infty \ .
\end{equation}
\item {\it Existence of a regular cyclic vector:}
\begin{equation}
\la{2.2}
\exists\,i \in K^0 \ \mbox{such that}\ m_2^0(i, \,\mbox{--}\,):
\, A \stackrel{\sim}{\to} K^0 \ \mbox{is an isomorphism of vector spaces.}
\end{equation}
\item {\it Weak associativity:}\quad For all $ a \in A $ and for all $ v \in K^1\,$ we have
\begin{equation}
\la{2.3}
m_3(v, x, a) = 0 \ ,
\end{equation}
\begin{equation}
\la{2.4}
m_3(v, a, y) = 0 \ ,
\end{equation}
\begin{equation}
\la{2.5}
m_3(v, y, x) \in k.i \ ,
\end{equation}
where $\, k.i \,$ denotes the subspace of $ K^0 $ spanned by the cyclic vector $ i \,$.
\end{itemize}
\end{definition}

\noindent
A few informal comments on these axioms may be relevant.

1. Since $ M $ is a $0$-complex, the quasi-isomorphism $\, r \,$ is  strict, and hence
$A$-linear. Moreover, since $\K$ has only two components, $\,r$ induces an isomorphism
of $A$-modules: $\, M \stackrel{\sim}{\to} H^0(\K) = \Ker(m_1) \,$, and the map
$ \, m_1: K^0 \to K^1  \,$ is surjective\footnote{The surjectivity of $ m_1 $ is also a
formal consequence of axiom \eqref{2.1} as the latter implies
$\, \dim\ \Coker(m_1) < \infty \,$ while $ A $ has no nontrivial finite-dimensional modules.}.
Now, $ M $ has only trivial ($A$-linear)
automorphisms, i.e. $\, \Aut_{A}(M) = k^{\times} \,$. Hence, being strict, the $\A$-morphism
$\, r \,$ is determined uniquely (up to a constant factor) by its target $ \K $.
Thus, we may (and often will) refer to $ \K \,$, rather than $ r \,$, as an $\A$-envelope of $M\,$.
See also Lemma~\ref{L3.1} below.

2. The axiom (\ref{2.2}) suggests to think of $ K^0 $ as a ``free module of rank $1$''
over $A\,$, though with $A$ acting non-associatively. Then, being a finite quotient of $ K^0 \,$,
$\,K^1 $ might be regarded as (a non-associative analogue of) a ``finite-dimensional cyclic
representation'' of $A\,$. Proposition~\ref{P1} below justifies in part this interpretation.

3. The axioms (\ref{2.3}) and (\ref{2.4}) together with structure relations (\ref{1.3})
imply
\begin{equation}
\la{2.3e}
m_3(v, \, k[x],\, A) \equiv 0 \ , \quad m_3(v, \, A, \, k[y]) \equiv 0 \ .
\end{equation}
These could be interpreted by saying that the elements of $ k[x] $ act
associatively on $\K$ when written ``on the left'', while the elements of $ k[y] $ act
associatively when written ``on the right'', i.~e.
$$
m_2(m_2(\,\mbox{--}\,,\, p),\, a) = m_2(\,\mbox{--}\,,\, p\,a) \ , \quad
m_2(m_2(\,\mbox{--}\,,\, a),\,q) = m_2(\,\mbox{--}\,,\, a\,q)
$$
for all $ p \in k[x]\,, q \in k[y] $ and $ a \in A\,$.

4. All the axioms above make sense in the commutative situation, and it is instructive
to see what happens if we replace the Weyl algebra in Definition~\ref{D1} by its
polynomial counterpart.
\begin{proposition}
\la{P1}
Suppose (for a moment) that $ A = k[x,y] \,$. If $\, \K = K^0 \oplus K^1
\in \Modi(A) \,$ satisfies \eqref{2.1}--\eqref{2.5} and $ m_1 $ is surjective
then $\, m_3 \equiv 0 \,$ on $\K\,$.
\end{proposition}
\begin{proof}
It suffices to show that $ m_3(v,y,x) = 0\,$ for all $\, v \in K^1$.
The vanishing of $ m_3 $ follows then routinely from \eqref{2.3e} and
commutativity of $A\,$. If $ K^1 = 0 $ there is nothing to prove.
So we may assume $ K^1 \not= 0\,$.
Then $ m_1(i) \not= 0 \,$ for $\, m_2^1(m_1(i), \,\mbox{--}\,) =
m_1 m_2^0(i, \,\mbox{--}\,): A \to K^1 \,$ is surjective by \eqref{2.2}.
Now, using the notation \eqref{2.7} -- \eqref{2.10} and arguing as in
Lemma~\ref{L3} below, we can compute $\, [\X, \Y] = \i\,\j \,$.
On the other hand, the set of vectors $ \{\, \Y^m\X^k(\i) \,\} $ spans $ K^1 $
and $\, \dim\, K^1 < \infty \,$. An elementary lemma from linear algebra
(see, e.g., \cite{N}, Lemma~2.9) implies then $ \j \equiv 0 $.
\end{proof}

\noindent
Thus, if $ A = k[x,y] \,$, an $\A$-module $ \K $ satisfying the axioms of Definition~\ref{D1}
can be identified with a usual complex of $A$-modules, $ K^0 $ being isomorphic to the
free module of rank $1$ and $ K^1 $ being a finite-dimensional cyclic representation of $ A $.
As mentioned in the Introduction, the latter corresponds
canonically to a point of the Hilbert scheme $\, \Hilb_{n}(\AA^2_k)\,$ with
$\, n = \dim\, K^1\,$.

Returning now to the Weyl algebra, we will see that the points of the Calogero-Moser varieties
$ \CC_n $ arise from $\A$-envelopes in a similar manner.

\subsection{The Calogero-Moser data}
\la{CMD}
Let $ \K $ be an $\A$-module satisfying the
axioms (\ref{2.1}) -- (\ref{2.5}). Denote by $\, X, Y \,$ (resp., $\, \X, \Y $) the action
of the canonical generators of $ A $ on $ K^0 $ (resp., $ K^1 $), i.~e.
\begin{equation}
\la{2.6}
X := m_2^0(\,\mbox{--}\,, x) \in \End_{k}(K^0)\ , \quad
Y := m_2^0(\,\mbox{--}\,, y) \in \End_{k}(K^0)\ ,
\end{equation}
\begin{equation}
\la{2.7}
\X := m_2^1(\,\mbox{--}\,, x) \in \End_{k}(K^1)\ , \quad
\Y := m_2^1(\,\mbox{--}\,, y) \in \End_{k}(K^1)\ .
\end{equation}
In view of (\ref{1.2}) we have
\begin{equation}
\la{2.8}
\X \, m_1 = m_1 \, X \ ,\quad \Y \, m_1 = m_1 \, Y \ .
\end{equation}
Now the axiom (\ref{2.5}) yields a $k$-linear functional $ \j $ on $ K^1 \,$ such that
\begin{equation}
\la{2.9}
m_3(v, y, x) = \j(v)\,i \quad \mbox{for all}\ v \in K^1\ .
\end{equation}
Combining $ \j $  and the cyclic vector $\, i \in K^0 \,$ (see (\ref{2.2})) with differential
on $\K $ we define
\begin{equation}
\la{2.10}
\i := m_1(i) \in K^1 \ , \quad j := \j \, m_1\, \in \Hom_{k}(K^0, k)\ .
\end{equation}
\begin{lemma}
\la{L3}
The data introduced above satisfy the equations:
\begin{equation}
\la{2.11}
X\,Y - Y \,X + \id_{K^0} = i\,j \ , \quad  \X\,\Y - \Y\,\X + \id_{K^1} = \i\, \j\ .
\end{equation}
\end{lemma}
\noindent
Indeed, in view of \eqref{2.8} and surjectivity of $ m_1 \,$, the second equation in
(\ref{2.11}) is a consequence of the first, while the first follows formally from \eqref{2.4} and
\eqref{2.5}:
\begin{eqnarray}
&& u = m_2^0(u, \, 1) = m_2^0(u, xy-yx) = m_2^0(u, xy) - m_2^0(u, yx)   \nonumber\\*[1ex]
&& \quad \! = m_2^0(m_2^0(u,x),y) - m_3(m_1(u), x, y) - m_2^0(m_2^0(u,y),x) + m_3(m_1(u), y, x)
\nonumber\\*[1ex]
&& \quad \! = Y\,X(u)  - X\,Y(u) + \j\,m_1(u)\,i  =  (Y\,X - X\,Y + i\,j)\,u \ ,
\quad \forall u \in K^0\ . \nonumber
\end{eqnarray}

Thus, given an $\A$-envelope $ \K $, the quadruple $ (\X, \Y, \i,\j) $ represents a point of the
Calogero-Moser variety $ \CC_n \,$, where $ n = \dim\,K^1 \,$. Conversely, given a quadruple
$ (\X, \Y, \i, \j) $ satisfying \eqref{2.11}, we will show now how to construct an associated
$\A$-envelope.

\subsection{From Calogero-Moser data to $\A$-envelopes}
\la{Recon}
Let $ R := k \langle x,y \rangle $ be the free algebra on two generators. Denote by
$\,\tau: R \to R \,$,$\,a \mapsto a^\tau $, the canonical anti-involution acting
identically on $x$ and $y$.
(Thus, $\, x^\tau = x \,$, $\, y^\tau = y \,$ and $\, (ab)^\tau = b^\tau a^\tau \,$,
$\,\forall\, a,b \in R\,$.) Given a quadruple $\, (\X, \Y, \i, \j) \,$ representing a point
of $ \CC_n $, we introduce the linear functional
\begin{equation}
\la{fun}
\varepsilon:\, R \to k\ , \quad  a(x,y) \mapsto \j \,a^\tau(\X, \Y) \,\i\ ,
\end{equation}
and define the right action of $ R $ on $ k^n $ by $\, a(x,y) \mapsto a^\tau(\X, \Y) \in
\End_k(k^n)\,$. By \cite{W}, Lemma~1.3, $\, k^n $ becomes then a cyclic (in fact, irreducible)
module over $ R $ with cyclic generator $ \i\,$. We denote this module by $ K^1 $ and write
$\, m: R \to K^1 \,$ for the $R$-module homomorphism sending $\, 1 \mapsto \i \,$. More
explicitly, we have $\, m: a(x,y) \mapsto a^\tau(\X, \Y)\,\i $ and hence the equality
$\, \varepsilon = \j\,m \,$.

Next, we form the following right ideal in the algebra $ R $
\begin{equation}
\la{id}
J := \sum_{a \in R} (a \, w + \varepsilon(a))\,R\ ,
\end{equation}
where $ w := xy-yx-1 \in R\,$, and let $\, K^0 := R/J\,$. Clearly, $ K^0 $ is
a cyclic right module over $R$ whose generator $\, [\,1\,]_J \,$ we denote by $ i\,$.

Note that both maps $\, \varepsilon \,$ and $\, m \,$ factor through the canonical projection
$\, R \onto R/J \,$, thus defining a linear functional $ j: K^0 \to k $ and an $R$-module
epimorphism $ m_1: K^0 \onto K^1 \,$ respectively. Indeed, since $\, \varepsilon = \j\,m \,$
it suffices to check that $ m $ vanishes on $ J $, and that is an easy consequence of our
definitions:
$$
m(a\, w + \varepsilon(a)) = \left(w^\tau(\X, \Y) \,a^\tau(\X, \Y) + \varepsilon(a)\right) \i
= -\i\,\j\, a^\tau(\X, \Y)\,\i + \i\,\varepsilon(a) \equiv 0 \ .
$$
Now, we have obviously $\, \i = m_1(i) \,$ and $\, j = \j\,m_1 \,$. Moreover, if we
let $X$ and $Y$ denote the endomorphisms of $ K^0 $ coming from the action of
$ x $ and $ y $ in $ R $ then
$$
(XY - YX + \id)[\,a\,]_J = [\,a(yx - xy + 1)\,]_J = [\,- a\,w\,]_J =
[\,\varepsilon(a)\,]_J  = j([\,a\,]_J)\,i \ ,
$$
and hence the relation $\, XY - YX + \id = i\,j\,$.

Summing up, we have constructed a complex of vector spaces
$$
\K :=[\,0 \to K^0 \stackrel{m_1}{\longrightarrow} K^1 \to 0 \,] \ ,
$$
together with linear data $ (X, Y, i, j) $ and $ (\X, \Y, \i, \j) $
satisfying \eqref{2.8}, \eqref{2.10} and \eqref{2.11}. Clearly, assigning $ (X, \X) $
and $ (Y, \Y) $ to the canonical generators of $ A $ does not make $ \K $ a
complex of $A$-modules. However, this {\it does} define an action of $ A $ on
$ \K\,$ ``up to homotopy''. More precisely, we have
\begin{lemma}
\la{L6}
The assignment $\, x \mapsto [(X, \X)] \,$ and $\, y \mapsto [(Y, \Y)] \,$
extends to a well-defined {\rm algebra} homomorphism
\begin{equation}
\la{homot}
\alpha: \, A \to \End_{\HH(k)}(\K)^{\mathsf{opp}}\ ,
\end{equation}
where $ \HH(k) $ denotes the homotopy category of $ \Com(k) $.
\end{lemma}
\begin{remark}
Given an algebra map (\ref{homot}), we say that $ A $ acts {\it homotopically} on the
complex $ \K $ and refer to $ (\K, \, \alpha) $  as a (right) {\it homotopy module} over
$A\,$ (cf. \cite{K2}).
\end{remark}
\begin{proof}
We need only to check that $\, [x,\,y] \,$ acts on $ \K $ by an endomorphism
homotopic to the identity map. This is an easy consequence of \eqref{2.11}.
Indeed, we have
$$
\alpha([x,\,y]) - \id_{\K}  = ([Y, \,X] - \id_{K^0},\, [\Y, \,\X] - \id_{K^1}) =
(- i\,j,\, -\i\,\j) \ ,
$$
so the required homotopy $\, h: K^1 \to K^0  \,$ satisfying 
$\, h \circ m_1 = - i\,j \,$ and $\, m_1 \circ h = -\i \,\j \,$
is given by $\, h = -i\,\j\,$. 
\end{proof}

Let $\,\pi: \End_{\Com(k)}(\K)^{\mathsf{opp}} \to \End_{\HH(k)}(\K)^{\mathsf{opp}} \,$ be the 
canonical projection assigning to an endomorphism of $ \K $ its homotopy  class.
Thus, $ \pi $ is an algebra map with $ \Ker(\pi) $ consisting of null-homotopic endomorphisms.
Now, to make a homotopy module $\, (\K, \, \alpha) \,$ a unital $\A$-module over
$ A $ it suffices to choose a linear lifting
\[
\begin{diagram}[small, tight]
  &                                &   \End_{\Com(k)}(\K)^{\mathsf{opp}} \\
  & \ruTo^{\varrho}                &   \dOnto_{\pi}       \\
A &    \rTo^{\alpha}               & \End_{\HH(k)}(\K)^{\mathsf{opp}}    \\
\end{diagram}
\]
such that
\begin{equation}
\la{lift}
\pi \circ \varrho = \alpha  \quad \mbox{and} \quad \varrho(1) = \id_{\K}\ .
\end{equation}
Indeed, given such a lifting, we can define
\begin{eqnarray}
\la{oper}
&& m_1: \K  \to \K\ ,\quad (u,v) \mapsto (0,\, m_1(u))\ ,\nonumber\\
&& m_2: \K \otimes A \to \K\ ,\quad (u,v) \otimes a \mapsto (\varrho^0(a)u, \,
\varrho^1(a)v)\ ,\\
&& m_3: \K\otimes A \otimes A \to \K\ , \quad (u,v)\otimes a \otimes b
\mapsto (-\omega^0(a,b)\,m_1^{-1}(v),\,0) \ ,\nonumber
\end{eqnarray}
where $\, \omega: A \otimes A \to \End_{\Com(k)}(\K)^{\mathsf{opp}} \,$ denotes
the ``curvature'' of the map $ \varrho $ which measures its deviation from being
a ring homomorphism (see \cite{Q1}):
\begin{equation}
\la{curv}
\omega(a,b) := \varrho(ab) - \varrho(b) \, \varrho(a) \ , \quad a, b \in A\ .
\end{equation}
Note, in view of (\ref{lift}), $\, \omega(a, b) \in \Ker(\pi)\,$ for all $\,a,b \in A\,$.
Hence $ \omega(a,b) $ is null-homotopic and therefore induces the zero map on cohomology
of $ \K \,$.
Since in our case $\, H^0(\K) = \Ker(m_1) \,$, we see that
$\, \omega^0(a,b): K^0 \to K^0 \,$ vanishes on $ \Ker(m_1) $ and thus
induces naturally a linear map $\, \omega^0(a,b)\,m_1^{-1}: K^1 \to K^0 \,$.
This justifies the definition of $ m_3 $ in (\ref{oper}).

It is now a trivial exercise to check that the maps (\ref{oper}) satisfy the first three
defining relations (\ref{1.1})--(\ref{1.3}) of $\A$-modules.
Since $ \K $ is a two-term complex with surjective $ m_1 \,$, Lemma~\ref{L1} guarantees then
that $ \K $ is a genuine $ \A$-module over $ A\,$.  Moreover, $ \K $ is unital due to the
last condition in (\ref{lift}). Thus, we need only to find a specific lifting that would
verify the axioms \eqref{2.2}--\eqref{2.5}.

There is an obvious choice for such a lifting: namely, we may define $\,\varrho \,$ by
\begin{equation}
\la{lift1}
\varrho(x^k y^m) := (Y^m X^k,\ \Y^m \X^k)\ , \quad \forall\, k,m \geq 0\ .
\end{equation}
Then, choosing the monomials $\,\{x^k y^m\}\,$ as a linear basis in $ A $, we have
$$
m_2^0(i,\,x^k y^m) = Y^m X^k(i) = [\,x^k y^m\,]_{J}\,\in K^0 \ ,
$$
where $ [\,x^k y^m\,]_{J} $ denotes the residue class of $\, x^k y^m  \in R \,$
modulo  $ J\,$. Such residue classes are all linearly independent and span
$ R/J $ as a vector space.  Hence, $\, m_2^0(i,\,\mbox{---}\,): A \to K^0 = R/J \,$
is a vector space isomorphism as required by \eqref{2.2}. The conditions (\ref{2.3})--(\ref{2.5})
are verified at once by computing the ``curvature'' of \eqref{lift1} and substituting the
result in \eqref{oper}: for example,
$$
\omega^0(y, x) = \varrho^0(yx) - \varrho^0(x) \, \varrho^0(y)  = \varrho^0(xy-1) - XY
= YX - XY - \id_{K^0}  = - i\,j\ ,
$$
and hence $ m_3(v, y, x) = \j(v)\,i $ for all $ v \in K^1\,$.

Thus, starting with Calogero-Moser data $ (\X, \Y, \i, \j) $, we have constructed an $\A$-module
$ \K $ that satisfies the axioms of Definition~\ref{D1}. It remains only to show that $ \K $ represents a rank $1$ torsion-free $A$-module in $\Di(A)$.
\begin{lemma}
\la{L3.1}
If $ \K \in \Modi(A) $ satisfies \eqref{2.1}--\eqref{2.5} then $ H^0(\K) $ is a
finitely generated rank $1$ torsion-free module over $A$.
\end{lemma}
\begin{proof} Fix some standard increasing filtration on $ A \,$, say
$\, A_n := \Span\{x^k y^m \,:\, m+k \leq n\}\,$, so that $\, \grd(A) := \oplus_{n\geq 0}
A_n /A_{n-1} \cong k[x,y]\,$. With isomorphism \eqref{2.2} we can transfer this
filtration on the complex $ \K \,$: more precisely, we set
$\, K^0_{n} := m_2^0(i,\, A_n) \,$ and $\, K^1_{n} := m_2^1(\i,\, A_n) \,$ for each $ n \geq 0\,$.
Now, using the relations
\eqref{2.11} it is easy to see that $ m_2(\K_n, A_m) \subseteq \K_{n+m} $ for all $ n,m \geq 0\,$.
Hence, the $\A$-structure on $ \K $ descends to the associated graded complex $\,
\grd(\K) := \oplus_{n\geq 0} \K_n/\K_{n-1} \,$ making it an $ \A$-module over $ \grd(A)\,$.
Relative to $ \grd(A)\,$, this module satisfies the same axioms \eqref{2.1}--\eqref{2.5}
as $ \K $, and hence by Proposition~\ref{P1}, it must be a genuine complex of $ \grd(A)$-modules.
In particular, we have $\, \grd(K^0) \cong \grd(A) \,$ (as $ \grd(A)$-modules).
Putting now on $ H^0(\K) = \Ker(m_1) \subseteq K^0 $ the induced filtration and passing to the associated graded level we see that $ \grd\,H^0(\K) $
is a f.~g. rank $1$ torsion-free module over $\grd(A)$ (as it canonically embeds in $ \grd(K^0) $).
By standard filtration arguments all the above properties lift to $ H^0(\K) \,$. Hence
$ H^0(\K) $ is a f.~g. rank $1$ torsion-free module over $A$.
\end{proof}

\section{Envelopes vs. Resolutions}
\la{MRes}
In this section we show how to construct some explicit representatives of (the isomorphism
class of) a module $ M $ from its $\A$-envelope $\, M\stackrel{r}{\longrightarrow} \K \,$.
The key idea is to relate $ \K $ to a {\it minimal} injective resolution of $ M \,$.

Thus, let $\, e: M \to \E \,$ be a minimal injective resolution of
$ M \,$ in $ \Mod(A)\,$. This has length one, i.~e. $\, \E =
[\,0 \to E^0 \stackrel{\mu_1}{\longrightarrow}  E^1 \to 0\,]\,$, and is uniquely determined
(by $ M $)  up to isomorphism in $ \Com(A)\,$.
When regarded as an object in $ \Modi(A) $, $\, \E \,$ represents the same quasi-isomorphism
class as
$ \K \,$. It is therefore natural to find a quasi-isomorphism that ``embeds''
$ \K $ in $ \E\,$. Indeed, if $ \K $ {\it were} a genuine
complex of $ A$-modules, such an embedding would always exist in $ \Com(A) $
and would be unique and canonical by injectivity of $ \E\,$. In our situation, however,
no  {\it strict}\, quasi-isomorphism in $ \Modi(A) \,$ maps
$\, r \,$ to $\, e \,$ (unless $ M $ is free). Instead, we will construct two
``partially strict'' quasi-isomorphisms $\, g_x: \K \to \E \,$ and  $\, g_y: \K \to \E \,$,
the first being linear with respect to the action of $ k[x] $ and the second with respect
to the action of $ k[y]\,$. As we will see, such maps are unique and defined canonically
(depending only on the choice of generators $x$ and $y$ of the algebra $A$).
What seems remarkable is that both $ g_x $ and $ g_y $ can be expressed explicitly in
terms of the Calogero-Moser data. Identifying then $ E^0 $ with $ Q $
(the quotient field of $A$) and restricting our maps to the cohomology of $ \K $ we
will get two distinguished representatives of $ M $ as fractional ideals in $Q$.

Before stating our main theorem we notice that any $\A$-morphism $\, g: \K\to \E \,$
has at most two nonzero components: with a slight abuse of notation, we will write these
in the form
\begin{eqnarray}
&& g_1: \K \to \E\ , \quad (u,\,v) \mapsto (g_1(u), \,\g_1(v))\ ,\nonumber\\
&& g_2: \K \otimes A \to \E\ ,\quad (u,\,v) \otimes a \mapsto (g_2(v,a),\, 0)\ .\nonumber
\end{eqnarray}
\begin{theorem}
\la{T2}
Let $\,r: M \to \K $ be an $ \A$-envelope of $ M \,$,
and let $\, e: M \to \E \,$ be a minimal injective resolution of $M$ in $ \Mod(A)\,$.

$ (a) $\  There is a unique pair $\, (g_x,\, g_y)\,$ of $\, \A$-quasi-isomorphisms
making the diagram
\begin{equation}
\la{D3}
\begin{diagram}[small, tight]
   &           & M &                    &    \\
   & \ldTo^{r} &   & \rdTo^{e}     \\
\K & &  \pile{\rTo^{g_x}\\ \rTo_{g_y}} & & \E \\
\end{diagram}
\end{equation}
commutative in $\, \Modi(A) \,$ and satisfying the conditions
\begin{equation}
\la{2.16}
(g_x)_2\,(v, x) = 0 \quad \mbox{and} \quad (g_y)_2\,(v, y) = 0 \ ,\quad  \forall\, v \in K^1 \ .
\end{equation}

$(b)$\ If we choose $\,\{ Y^m X^k i\}\, $ as basis in $ K^0 \,$
and $\,\{x^k y^m\} \, $ as basis in $ A $ then $\, g_x \,$ and $\, g_y \,$ are given by
\begin{eqnarray}
&&(g_x)_1\,(\,Y^m X^k i\,) =
i_x \cdot \left(\,x^k y^m + \Delta_{x}^{k m}(\i) \,\right)\ , \quad
(g_x)_2\,(v,\,x^k y^m) =  i_x\cdot \Delta_{x}^{k m}(v)\ , \la{2.21x}\\*[1ex]
&&(g_y)_1\,(\,Y^m X^k i\,) = i_y \cdot \left(\,x^k y^m + \Delta_{y}^{k m}(\i)\,\right)\ , \quad
 (g_y)_2\,(v,\,x^k y^m) =  i_y\cdot \Delta_{y}^{k m}(v) \ ,\la{2.21y}
\end{eqnarray}
where $\, i_x := (g_x)_1\,(\,i\,) \,$ and $\, i_y := (g_y)_1\,(\,i\,)\,$ in $ E^0 $, and
\begin{eqnarray}
&& \Delta_{x}^{k m}(v) := - \, \j(\X - x)^{-1} (\Y-y)^{-1} (\Y^m - y^m)\,\X^k v \ ,
\la{2.22x}\\*[1ex]
&& \Delta_{y}^{k m}(v) :=  \j(\Y - y)^{-1} (\X-x)^{-1} (\X^k - x^k)\,y^m v \ .\la{2.22y}
\end{eqnarray}
Moreover, we have
\begin{equation}
\la{ii}
i_x = i_y \cdot \kappa\ ,
\end{equation}
where
\begin{equation}
\la{kappa}
\kappa := 1 - \j\,(\Y - y)^{-1} (\X-x)^{-1}\i \ \in Q \ .
\end{equation}
\end{theorem}

\vspace{1ex}

\noindent
Part $(b)$ needs perhaps some explanations.

1.\ The set $\,\{Y^m X^k i\,:\,m,k \geq 0\}\, $ is indeed a linear basis in $ K^0 $
because it is the image of the linear basis $\,\{x^k y^m \,:\, m,k \geq 0 \} \, $ of $ A $
under the isomorphism (\ref{2.2}).

2.\ The formulas (\ref{2.22x}) and (\ref{2.22y}) define the
maps $\, \Delta_{x,y}^{k m}:\, K^1 \to Q \,$ for $ m,k \geq 0\,$, which
could be written more accurately as follows
\begin{eqnarray}
&&\Delta_{x}^{k m}(v) := - \det(\X - x\,\id)^{-1}
(\j \otimes 1)\,[\,(\X - x \,\id)^{*}
\sum_{l=1}^{m}\, \Y^{m-l}\X^k(v) \otimes y^{l-1}\,] \ ,\nonumber\\
&&
\Delta_{y}^{k m}(v) := \det(\Y - y\,\id)^{-1}
(\j \otimes 1)\,[\,(\Y - y \,\id)^{*}
\sum_{l=1}^{k}\, \X^{k-l}(v) \otimes x^{l-1}y^m\,] \ ,\nonumber
\end{eqnarray}
where $\, \id := \id_{K^1} \,$, $\, (\X - x\,\id)^{*} \in \End_{k}(K^1) \otimes A \,$
denotes the classical adjoint of the matrix $\, \X - x\,\id \,$ and
$ \j \otimes 1\,: K^1 \otimes A \to A $ is defined naturally by
$\, v \otimes a \mapsto \j(v)\,a \,$.

3.\ The dot in the right hand sides of (\ref{2.21x}) and (\ref{2.21y}) denotes the (right)
action of $ A $ on $ \E \,$. Even though $\, \Delta_{x,y}^{k m}(v) \in Q \,$, these
formulas make sense since both components of $ \E $ are injective (and hence
divisible) modules over $ A\,$.\\

Now we proceed to the proof of Theorem~\ref{T2}. We will describe in detail only
the map $ g_x $ writing it simply as $ g\,$. Repeating a similar construction for
$ g_y $ is a (trivial) exercise which we will leave to the reader.

First, observe that (\ref{2.3e}) implies $\, m_3(v,\,k[x],\,k[x]) \equiv 0 \,$,
and thus allows us to treat
$\, \K \,$ as a usual complex of $k[x]$-modules (via the embedding $ k[x] \hookrightarrow A $).
Being strict, the quasi-isomorphism
$\, r: M \to \K \,$ is $k[x]$-linear, and hence can also be regarded as a quasi-isomorphism
in $ \Mod(k[x])\,$. Now, since $ A $ is projective (in fact, free) as $k[x]$-module, every
injective over $A$ is automatically injective over $ k[x] $ (see \cite{CE}, Proposition~6.2a, p.~31).
Hence, $ e: M \to \E $ extends to a $k[x]$-linear
morphism $\, g_1: \K \to \E \,$  such that the diagram
\begin{equation}
\la{D4}
\begin{diagram}[small, tight]
0 & \rTo &  M   & \rTo^{r}           & K^0         & \rTo^{m_1} & K^1 & \rTo   & 0 \\
  &      &  \dEq &                    & \dTo^{g_1}&            & \dTo^{\g_1} &   \\
0 & \rTo &  M    & \rTo^{e} & E^0         & \rTo^{\mu_1} & E^1 & \rTo   & 0 \\
\end{diagram}
\end{equation}
\noindent
commutes in $ \Com(k[x])\,$. We claim that such an extension is unique.
Indeed, if $\, g_{1}': K^0 \to E^0 \,$ is another map in $ \Mod(k[x]) $
satisfying $\, g_1 \circ r = g_1' \circ r  = e \,$,
then $\,  g_1' - g_1  \equiv 0 \,$ on $ \Ker(m_1)\,$ by exactness of the
first row of (\ref{D4}). So the difference
$\, d := g_1' - g_1 \,$ induces a $k[x]$-linear map $\,
\bar{d}: K^1 \to E^0 \,$. Since $\, \dim_k\,K^1 < \infty \,$,
$\, K^1 $ is torsion over $ k[x]\,$, while $ E^0 $ is obviously torsion-free.
Hence, $\, \bar{d} = 0 \,$ and therefore $\, g_1' = g_1 \,$.
This implies, of course, that $\, g_1'  = g_1 \,$ as morphisms in $\Com(k[x])$.

\begin{lemma}
\la{L4}
The map $\, g_1: \K \to \E \,$ extends to a unique quasi-isomomorphism
of $\A$-modules over $ A\,$.
\end{lemma}
\begin{proof}
According to Lemma~\ref{L2}, it suffices to show the existence of a map
$$
g_2: K^1 \otimes A \to E^0
$$
satisfying the conditions (cf. (\ref{1.6e}) and (\ref{1.61e}))
\begin{equation}
\label{2.12}
g_1(m_2^0(u,a)) - g_1(u) \cdot a = g_2 (m_1(u),a) \ ,
\end{equation}
\begin{equation}
\label{2.13}
\g_1(m_2^1(v,a)) - \g_1(v) \cdot a  = \mu_1(g_2(v,a)) \ .
\end{equation}
Since $ m_1 $ is surjective, $ (\ref{2.13}) $ is a consequence of (\ref{2.12}) and
commutativity of the diagram (\ref{D4}). On the other hand, to satisfy (\ref{2.12})
we need only to show
\begin{equation}
\label{2.14}
g_1(m_2^0(u,a)) = g_1(u) \cdot a \quad \mbox{for all}\ u \in \Ker(m_1)\ ,
\end{equation}
and this again follows easily from the diagram (\ref{D4}). Indeed, since the first row is
exact, we have $\,u = r(m)\,$ for some $ m \in M\,$ whenever $\, u \in \Ker(m_1) \,$, and
in that case
$\,
g_1(u) \cdot a = g_1\,r(m) \cdot a = e(m)\cdot a =
e(m \cdot a) = g_1\,r(m \cdot a) = g_1(m_2^0(u,a))\,$.
Thus, we can simply define $ g_2 $  by the formula
\begin{equation}
\label{2.15}
g_2(v,a) := g_1(m_2^0(m_1^{-1}(v),a)) - g_1(m_1^{-1}(v)) \cdot a\ .
\end{equation}
which makes sense due to (\ref{2.14}). The uniqueness of $ g_2 $ is obvious.
\end{proof}

Clearly, the $\A$-morphism given by Lemma~\ref{L4} satisfies the conditions on
$ g_x $ of Theorem~\ref{T2}$(a)$: in fact, $\, g \,$
being $ k[x]$-linear means
\begin{equation}
\la{2.16e}
g_2(v, \, x^k) = 0 \quad \forall\, v \in K^1\ , \ \forall\, k \geq 0 \ .
\end{equation}

On the other hand, if an $ \A$-morphism  $ g: \K \to \E $ satisfies
$\, g_2(v, x) = 0 \,$,$\,\forall\, v \in K^1 \,$, then (\ref{2.16e})
holds automatically. This is immediate by induction from (\ref{2.12}) and
the axiom (\ref{2.3}). Thus, the uniqueness of $ g_x $ follows again from
Lemma~\ref{L4}. This finishes the proof of Part $(a)$ of the Theorem.\\

To prove Part $(b)$ we start with the identity $\, g_3 = 0 \,$ which
holds automatically once the existence of the $\A$-morphism $ g\,$ is
established. As in Lemma~\ref{L2}, we will regard this identity as an equation
on $\,g_1\,$ and $\, g_2 \,$.
Taking into account that $\, m_3 \equiv 0 \,$ on $ \E \,$, we can write it in the
form (cf. ({\ref{1.11e})):
\begin{equation}
\la{2.17}
g_2(v,a) \cdot b = g_{2}(v,ab) - g_2(m_{2}^1(v,a),b) + g_1(m_3(v,a,b))\ ,
\end{equation}
where $\, v \in K^1\,$ and $\, a,b \in A\,$. Now, it turns out that (\ref{2.17})
can be solved easily, by elementary algebraic manipulations.

First, letting $\, a = x^k \,,\, b = y^m \,$ in (\ref{2.17}) and using  (\ref{2.3e}) and
(\ref{2.16e}), we get
\begin{equation}
\la{2.171}
g_2(v,\,x^k y^m) = g_{2}(\X^k(v), \, y^m) \quad \mbox{for all}\ m,k \geq 0\ .
\end{equation}
Next, we substitute $\, a = y \,$ and $ \, b = x \,$ in (\ref{2.17}) and
use (\ref{2.9}). Since $\, g_{2}(v, yx) = g_{2}(v, xy) = g_{2}(\X(v), y) \,$
by (\ref{2.171}), we have
$$
g_2(v,y) \cdot x - g_{2}(\X(v), y) =  \j(v)\,i_x \ ,
$$
where $\, i_x := g_1(i) \in E^0 \,$. This equation has a unique solution (otherwise the
difference of two solutions would provide a nontrivial $k[x]$-linear map: $ K^1 \to E^0 $
which is impossible), and it is easy to see that that solution is given by
\begin{equation}
\la{2.19}
g_2(v,y) = - \, i_x \cdot \j\,[\,(\X - x)^{-1} \,v\,]\ .
\end{equation}
Finally, with $\, a = y^{m-1} \,$ and $\, b = y \,$ (\ref{2.17}) becomes
the recurrence relation
$$
g_2(v,\,y^m) = g_2(v,\,y^{m-1}) \cdot y + g_2(\Y^{m-1}(v),\,y)\ ,
$$
which sums up easily
\begin{equation}
\la{2.172}
g_2(v,\,y^m) = \sum_{l=1}^{m}\, g_2(\Y^{m-l}(v),\,y)\cdot y^{l-1} \ ,\quad m \geq 1\ .
\end{equation}
Combining (\ref{2.171}), (\ref{2.19}) and (\ref{2.172}) together, we find
\begin{equation}
\la{2.18}
g_2(v,\,x^k y^m) = -\, i_x \cdot \j(\X - x)^{-1} (\Y-y)^{-1} (\Y^m - y^m)\,\X^k v\ ,
\end{equation}
which is exactly the second formula of \eqref{2.21x}; the first one follows now from \eqref{2.12}:
\begin{eqnarray}
\lefteqn{g_1(Y^m X^k i) = g_1(m_2^0(i,\, x^k y^m)) =
g_1(i)\cdot x^k y^m + g_2(\i,\,x^k y^m)  } \nonumber \\*[1ex]
&& \ \ \ \ \ \ \ \ \ \ = i_x \cdot
\left(\,x^k y^m - \j(\X - x)^{-1} (\Y-y)^{-1} (\Y^m - y^m)\,\X^k\,\i\,\right)\ .
\nonumber
\end{eqnarray}
A similar calculation (with roles of $ x $ and $ y $ interchanged) leads to
formulas \eqref{2.21y}.

The relation \eqref{ii} can be deduced from \eqref{2.21x} and \eqref{2.21y} as follows.
First, we observe that  $\, (g_x)_1 = (g_y)_1 \,$ on $\, \im(r)\,$, which
is immediate in view of commutativity of the diagram \eqref{D3}.
Now, by the Hamilton-Cayley theorem, the polynomial
$\, p(x) := \det(\X - x) \,$ acts trivially on $ K^1 \,$,
i.~e. $\, m_2^1(v,\, p(x)) = p(\X)v = 0 \,$ for all $\, v \in K^1 \,$.
Hence
\begin{equation}
\la{van}
m_1(p(X)i) = m_1 m_2^0(i,\, p(x)) = m_2^1(\i,\, p(x)) = 0\ ,
\end{equation}
and therefore $\, p(X) i \in \Ker(m_1) = \im(r)\,$. Thus, we have
\begin{equation}
\la{neweq}
(g_x)_1\left(p(X) i\right) = (g_y)_1\left(p(X) i\right)\ .
\end{equation}
By \eqref{2.21x}, the left hand side of \eqref{neweq} is $\, i_x \cdot p(x) \,$.
On the other hand, \eqref{2.21y} together with the identity $\, p(\X)v = 0 \,$ yields
\begin{eqnarray}
(g_y)_1\left(\,p(X) i\,\right) &=&
i_y \cdot \left[p(x) + \j(\Y - y)^{-1}(\X -x)^{-1}(\,p(\X) - p(x)\,)\,\i\,\right]  \nonumber \\*[1ex]
&=& i_y \cdot \left[1 - \j(\Y - y)^{-1}(\X -x)^{-1}\i\,\right]\,p(x) \ .\nonumber
\end{eqnarray}
Now, since $ E^0 $ is a torsion-free $A$-module, the equation \eqref{neweq} implies \eqref{ii}.
This finishes the proof of Theorem~\ref{T2}.\\

As an application of Theorem~\ref{T2}, we can describe the cohomology
of an $\A$-envelope in terms of its Calogero-Moser data.
\begin{corollary}
\la{C1}
Let $\, \K \in \Modi(A) \,$ be an $\A$-envelope of $ M $,
and let $\, (\X, \Y, \i, \j) \,$ be the Calogero-Moser
data associated with $ \K $. Then, $\, M \,$ is isomorphic to each of the
following (fractional) ideals
\begin{eqnarray}
&& M_x := \det(\X - x)\,A + \kappa^{-1} \,\det(\Y - y)\,A \ , \la{2.Mx} \\*[1ex]
&& M_y := \det(\Y - y)\,A + \kappa \,\det(\X - x)\,A\ , \la{2.My}
\end{eqnarray}
where $\, \kappa \in Q\, $ is given by formula \eqref{kappa}.
\end{corollary}
\begin{remark}
\la{R4}
It is easy to see that $ M_x $ and $ M_y $ are the ``distinguished representatives'' of
(the isomorphism class of) $ M $ in the sense of \cite{BW2} (see {\it loc. cit.}, Section~5.1).
\end{remark}
\begin{proof}
Recall that an injective envelope of a rank one torsion-free
module over a Noetherian domain is isomorphic to its quotient field
(see, e.g. \cite{B}, Exemple~1, p.~20). Thus, if $\, \E \,$ is a minimal injective
resolution of $ M $, and $\, g: \K \to \E $ is one of the maps
constructed in Theorem~\ref{T2}, there is a (unique) $A$-module
isomorphism $\, E^0 \stackrel{\sim}{\to} Q \,$ sending
$\,  g_1(i) \in E^0 \,$ to $\, 1 \in Q \,$.  Using this isomorphism  we
can identify $ E^0 $ with $ Q $ and compute the image of $ H^0(\K) = \Ker(m_1) $
under $ g_1 $ with the help of Theorem~\ref{T2}.
As a result, for $ g = g_x $ we will get the ideal $ M_x \,$, and
for $ g = g_y $ the ideal $ M_y \,$. We will consider only $ g = g_x $
leaving, as usual, $\, g_y $ to the reader.

Let $\, p(x) := \det(\X - x) \,$ and $\, q(y) := \det(\Y - y) \,$.
Then,  $\, p(X) i \in \Ker(m_1) \,$  by \eqref{van}, and similarly $\, q(Y) i \in \Ker(m_1) \,$.
We claim that these elements generate $ \Ker(m_1) $ as $A$-module.
Indeed, the submodule $\, m_2^0(p(X) i,\, A) + m_2^0(q(Y) i,\, A) \subseteq \Ker(m_1) \,$
has finite codimension in $ K^0 $, and hence {\it a fortiori}\, in $ \, \Ker(m_1) \,$.
But the latter is a genuine $A$-module and therefore cannot have proper submodules
of finite codimension. It follows that
\begin{equation}
\la{Ker}
\Ker(m_1) = m_2^0(p(X) i,\, A) + m_2^0(q(Y) i,\, A)\ .
\end{equation}
Thus, it suffices to compute the images of $ p(X)i $ and $ q(Y)i $
under $ g_1\,$. Such a computation has already been done in the proof of Theorem~\ref{T2}:
the image of $ p(X)i $ is given by $\, i_x \cdot p(x) \,$, and
$\, g_1\left(q(Y)i\right) = i_x \cdot \left[1 + \j\,(\X - x)^{-1} (\Y-y)^{-1}\i\right] \, q(y) \,$.
Now, if we identify $\, E^0 \cong Q \,$ (letting $\, i_x \mapsto 1 \,$) then by \eqref{Ker}
\begin{equation}
\la{2.27}
g_1 (\,\Ker(m_1)\,) =  p(x)\,A + \chi\, q(y)\,A \ ,
\end{equation}
where $\, \chi := 1 + \j\,(\X - x)^{-1} (\Y-y)^{-1}\i \, \in Q\,$. Using  \eqref{2.11},
is easy to check that $\, \chi\,\kappa = 1  \,$ in $ Q \,$, so the right hand side of
\eqref{2.27} is precisely $ M_x \,$.
\end{proof}

\section{Uniqueness}
\la{Un}
The aim of this section is to prove the uniqueness of $\A$-envelopes. As we will see,
the latter should be understood in the strong sense: to wit, the $\A$-envelopes are defined uniquely
up to unique strict isomorphism. The key result here (Theorem~\ref{T3}) establishes an
equivalence between different types of isomorphisms of $\A$-envelopes, and it is perhaps
the most important consequence of our axiomatics.

Before stating this theorem, we introduce some numerical invariants to distinguish
between different $ \A$-envelopes. Specifically, keeping the notation of Section~\ref{CMD}
we associate to an $\A$-module $ \K $ the linear form
\begin{equation}
\la{invar}
\lambda: A \to k\ , \quad \lambda(a) := j\, m_2^0(i,\,a) = \j\, m_2^1(\i,\,a) \ .
\end{equation}
Equivalently, $\, \lambda \,$ can be defined by its values on the basis of monomials in $A$:
\begin{equation}
\la{numb}
\lambda_{lk} := \lambda(x^k y^l) = j \,Y^l X^k \,i = \j \,\Y^l\X^k \,\i \ ,
\end{equation}
and thus is determined by the double-indexed sequence of scalars
$\,\{\lambda_{lk} \,: \,k,l \geq 0\}$.

\begin{theorem}
\la{T3}
Let $ \K $ and $ \KK $ be two $ \A$-modules satisfying \eqref{2.1}--\eqref{2.5}.
Then the following are equivalent:

$(a)\ $ $ \K $ and $ \KK $ are {\rm strictly} isomorphic,

$(b)\ $ $ \K $ and $ \KK $ are isomorphic,

$(c)\ $ $ \K $ and $ \KK $ are quasi-isomorphic,

$(d)\ $ $ \K $ and $ \KK $ determine the same functionals \eqref{invar}, i.~e.
$\, \lambda = \tilde{\lambda} \,$.
\end{theorem}
\begin{proof}
The implications $\, (a) \Rightarrow (b) \Rightarrow (c) \,$ are trivial.
It suffices only to show that $ \, (c) \Rightarrow (d) \, $ and $ \, (d) \Rightarrow (a) \,$.

If $ \K $ satisfies \eqref{2.1}--\eqref{2.5} then $ H^0(\K) $ is a f.~g. rank $1$
torsion-free $A$-module (see Lemma~\ref{L3.1}). By Corollary~\ref{C1},
$ H^0(\K) $ is then isomorphic to the fractional ideals $M_x$ and $M_y$ which are
related by $\, M_y = \kappa\,M_x \,$ (see \eqref{2.Mx}, \eqref{2.My}).
The element $ \kappa \in Q $ is given by the formula \eqref{kappa}. Developing the right
hand side of (\ref{kappa}) into the formal series:
\begin{equation}
\la{series}
1 - \j\,(\Y - y)^{-1} (\X-x)^{-1}\,\i = 1 - \sum_{l, k \geq 0}\,
\left(\,\j \,\Y^l\X^k \,\i\,\right) \, y^{-l-1} x^{-k-1}
\end{equation}
we notice that the coefficients of \eqref{series} are precisely the numbers (\ref{numb}).
Now, it is easy to see (cf. \cite{BW2}, Lemma~5.1) that $ \kappa $ is uniquely determined,
up to a {\it constant} factor, by the isomorphism class of $ H^0(\K)$.
Hence, if $ \K $ are $ \KK $ are quasi-isomorphic $ \A$-modules,
we have $\, \tilde{\kappa} = c \cdot \kappa \,$ for some $ c \in k\,$. Comparing the
coefficients of \eqref{series} yields at once $\, c = 1 \,$ and
$\, \tilde{\lambda}_{lk} = \lambda_{lk} \,$ for all $ l,k \geq 0\,$.
Thus, we conclude $\,(c) \Rightarrow (d) \,$.

Now, assuming $ (d) $ we construct a strict isomorphism $\, f: \K \to \KK \,$.
By definition, $\, f $ is given by two components: $\, f^0: K^0 \to \tilde{K}^0 \,$
and $\, f^1: K^1 \to \tilde{K}^1 \,$. In view of (\ref{2.2}), there is an obvious
candidate for the first one: namely, we can define $\, f^0\,$ by the commutative diagram
\begin{equation}
\la{3.1}
\begin{diagram}[small, tight]
   &           & A &                    &    \\
   & \ldTo^{\varphi} &   & \rdTo^{\tilde{\varphi}}     \\
K^0 & &  \rTo^{f^0} & & \tilde{K}^0 \\
\end{diagram}
\end{equation}
where $\, \varphi := m_2^0(i, \, \mbox{--}\,) \,$ and
$\, \tilde{\varphi} := \tilde{m}_2^0(\tilde{i}, \, \mbox{--}\,) \,$.
Then $\, f^0 \, $ is an isomorphism of vector spaces.
We need only to show that $\,(i)\ f^0 $ commutes with the action of $\, A \,$ and
$\,(ii)\ f^0(\Ker(m_1)) = \Ker(\tilde{m}_1)\,$. The latter condition allows one to
define an isomorphism $\, f^1: K^1 \to \tilde{K}^1 \,$ in the natural way
(so that $\, f^1\,m_1 = \tilde{m}_1 \, f^0 \,$) while the former guarantees
that $\, (f^0, f^1) \,$ is a strict map of $\A$-modules.

Now, $\,(ii)\,$ follows easily from $\,(i)\,$. To see this, let
$\, N := f^0(\Ker(m_1))\,$ and $\, \tilde{N} := \Ker(\tilde{m}_1)\,$
in $\, \tilde{K}^0 \,$. If $\, (i) \,$ holds then $\, N \,$ is closed under
the action of $\, \tilde{X} \,$ and $\, \tilde{Y}\,$, and so is obviously $ \tilde{N}\,$.
Moreover, both on $\, N \,$ and $\, \tilde{N} \,$, and therefore on
their sum $\, N + \tilde{N} \,$ the commutator
$\, \tilde{Y} \tilde{X} - \tilde{X} \tilde{Y} \,$ acts as identity.
Since
$$
\dim\, (N + \tilde{N})/\tilde{N} \leq \dim\, \tilde{K^0}/\tilde{N} =
\dim\, \tilde{K}^1 < \infty
$$
we have at once $\, N + \tilde{N} = \tilde{N} \,$, and therefore $\, N \subseteq \tilde{N} \,$.
On the other hand, using \eqref{2.11} we find
$$
\lambda(1) = \j(\i) = \tr(\i\,\j) = \tr\,\id_{K^1} = \dim\, K^1 \ .
$$
Hence, if $ (d) $ holds, then $\, \dim\, K^1 = \dim\,\tilde{K}^1 \,$, and
therefore $\,\dim\, \tilde{K}^0/N =  \dim\, \tilde{K}^0/\tilde{N}\,$. It follows that
$\, N = \tilde{N} \,$.

Thus, it remains to prove
\begin{equation}
\la{3.2}
f^0(m_2^0(u, \, b)) = \tilde{m}_2^0(f^0(u),\, b)\ ,\quad \forall\,u \in K^0\ , \quad
\forall\, b \in A \ .
\end{equation}
Note, if (\ref{3.2}) holds for $ b = x $ and $ b = y $ then, in view of (\ref{2.3e}),
it holds (by induction) for any powers $ x^k $ and $ y^m $, and more generally for
any element $\, b = x^k y^m \in A\,$. Now, for $\, b = y \,$, the equation (\ref{3.2}) is
immediate since both isomorphisms $ \varphi $ and $ \tilde{\varphi} $ in \eqref{3.1} are
$ k[y]$-linear (again due to \eqref{2.3e}).

Thus, it suffices to check (\ref{3.2}) only for $\, b = x \,$. To this end we fix the linear
basis $ \, \{Y^m X^k i\} \,$ in $ K^0 $ as in Theorem~\ref{T2},
and verify (\ref{3.2}) for each $\, u = Y^m X^k i \,$. First, observe
$$
X Y^m  - Y^m X + m Y^{m-1} =  \sum_{l=0}^{m-1} \, (Y^{m-l-1} i)\,j\, Y^l \ , \quad
\forall \, m \geq 0\ ,
$$
which follows easily by induction from (\ref{2.11}). Hence
\begin{eqnarray}
\la{3.3}
\lefteqn{m_2^0(Y^m X^k i, \, x) = X Y^m X^k i = }  \\
&& Y^m X^{k+1} i - m Y^{m-1} X^k i + \sum_{l=0}^{m-1} \, (j\, Y^l X^k \,i) \, Y^{m-l-1} i \ .
\nonumber
\end{eqnarray}
On the other hand,
\begin{equation}
\la{3.1e}
f^0(Y^m X^k i) = f^0(m_2^0(i,\, x^k y^m)) = \tilde{m}_2^0(\tilde{i}, \, x^k y^m) =
\tilde{Y}^m \tilde{X}^k \tilde{i} \ .
\end{equation}
Applying now $ f^0 $ to (\ref{3.3}) and using (\ref{3.1e}), we see that
$$
f^0(m_2^0(Y^m X^k i, \, x)) =
\tilde{m}_2^0(\tilde{Y}^m \tilde{X}^k \tilde{i}, \, x) =
\tilde{m}_2^0(f^0(Y^m X^k i), \, x)
$$
holds for all $\, k, m \geq 0 \,$ if and only if
$$
j\, Y^l X^k \,i =  \tilde{j}\, \tilde{Y}^l \tilde{X}^k \,\tilde{i}\ ,\quad
\forall \ l, k \geq 0 \ .
$$
In view of \eqref{numb} the latter conditions are equivalent to $ (d) $.
Thus, if $ (d) $ holds, the map $\, f^0 \,$ is $ A$-linear and induces a strict
isomorphism $ f: \K \to \KK\,$, implying $\,(a)$.
This finishes the proof of Theorem~\ref{T3}.
\end{proof}

The uniqueness of $\A$-envelopes is an easy consequence of Theorem~\ref{T3}.
\begin{corollary}
\la{U1}
Let $\,r: M \to \K \,$ and $\,\tilde{r}: M \to \KK \, $ be two $ \A$-envelopes
in the sense of Definition~\ref{D1}. Then there is a unique {\rm strict} isomorphism
of $ \A$-modules $\, f: \K \to \KK \,$  such that $\,\tilde{r} = f \circ r \,$ in $\,\Modi(A)\,$.
Thus, an $\A$-envelope of $M$ is determined uniquely up to (unique) strict isomorphism.
\end{corollary}
\begin{proof}
Once $ M $ is fixed, the quasi-isomorphism $\, r: M \to \K \,$ is uniquely determined by $ \K $
to a (nonzero) scalar factor (see remarks following Definition~\ref{D1}).
Hence, it suffices to have {\it any} strict isomorphism $\, f: \K \to \KK \,$ in
$ \Modi(A)\,$: multiplying $ f $ by an appropriate factor we can always achieve
$\,\tilde{r} = f \circ r \,$. Now, the existence of such an isomorphism is guaranteed by
implication $\, (c) \Rightarrow (a) \,$ of Theorem~\ref{T3}. The uniqueness is clear for
the difference of any two strict morphisms satisfying $\,\tilde{r} = f \circ r \,$
vanishes obviously on $ \im(r) $ and induces an $A$-linear map $\, K^1 \to \tilde{K}^0 \,$
which is also zero by torsion considerations.
\end{proof}

\section{Existence}
\la{Ex} In this section we give two different constructions of
$\A$-envelopes\footnote{Another, more geometric but less
elementary construction is given in the Appendix.}. The first
construction refines the elementary treatment of ideals in
\cite{BW2} and can be described in a nutshell as follows. Given a
rank one torsion-free $A$-module $M$, one cannot embed $ M $ in $
A $ as a submodule of finite codimension. However, as shown in
\cite{BW2}, there are two {\it different} embeddings $\, M \into A
\,$, one being a map of $ k[x]$-modules and the other of $
k[y]$-modules, which do have finite cokernels of the same
dimension. Using these embeddings, we construct two complexes of
vector spaces, each quasi-isomorphic to $ M $, but on which the
algebra $ A $ does not act in the usual (strict) sense. It turns
out, however, that these complexes can be ``glued'' together by a
natural linear isomorphism, and on the resulting complex one can
define a weak, {\it homotopic} action of $ A $. As in
Section~\ref{Recon}, this last action can then be enriched to a
full structure of $\A$-module giving an $\A$-envelope of $M$.

The second construction is also elementary, but it involves a priori no distinguished realization of $M$ in $A$.
Instead, we use an inductive procedure which somewhat resembles the construction
of minimal models (semi-free resolutions) in rational homotopy theory (see \cite{FHT}).
As in the case of minimal models, this procedure is far from being canonical --- it involves a lot of
choices --- but the uniqueness of Section \ref{Un} guarantees that the result is independent of any choices.

We start with formulating the main theorem of this section.
\begin{theorem}
\la{T4}
Every finitely generated rank $1$ torsion-free module over $ A $ has an
$\A$-envelope in $\, \Modi(A) \,$ satisfying the axioms of Definition~\ref{D1}.
\end{theorem}

\subsection{The first construction} As shown in \cite{BW2}, Section~5.1, the
isomorphism class of each rank $1$ torsion-free module in $\Mod(A)
$ contains a pair of {\it fractional} ideals $\, M_x \,$ and $\,
M_y \,$, which are uniquely characterized by a list of properties
and, in particular, such that $\, M_x  \subset k(x)[y] \,$ and $\,
M_y \subset k(y)[x] \,$. Despite being fractional, these ideals
can be embedded in $ A $ with the help of the following maps
\begin{eqnarray}
&& \rho_x :  k(x)[y] \to A \ ,\quad a(x)\, y^m  \mapsto  a(x)_{+}\, y^m\ ,
\la{rx} \\*[1ex]
&& \rho_y :  k(y)[x] \to A \ ,\quad a(y)\, x^m  \mapsto  a(y)_{+}\, x^m\ ,
\la{ry}
\end{eqnarray}
where ``$\,+\,$'' means taking the polynomial part of the corresponding rational function. As in \cite{BW2},
we write $\, r_x  \,$ and $\, r_y \,$ for the restrictions of these maps to $ M_x $
and $M_y$ respectively and denote by $\, V_x := A/r_x(M_x) \,$ and $\, V_y := A/r_y(M_y) \,$ the
corresponding cokernels. In this way we get two complexes of vector spaces
\begin{equation}
\la{Kxy}
\K_x := \left[\, 0 \to A \to V_x \to 0 \,\right]\quad \mbox{and} \quad
\K_y := \left[\, 0 \to A \to V_y \to 0 \,\right]\ ,
\end{equation}
together with quasi-isomorphisms $\, r_x: M_x \to \K_x \,$ and $\, r_y: M_y \to \K_y \,$.
By definition, $\, \rho_x \,$ is $k[y]$-linear and $\, \rho_y \,$ is $k[x]$-linear
with respect to the natural (right) multiplication-actions.
Hence, $ \K_x $  can be viewed as a complex of $ k[y]$-modules and
$ \K_y $ as a complex of $k[x]$-modules. Note, however, that neither on
$ \K_x $ nor on  $ \K_y $ the {\it full} algebra $A$ acts.

Now, since $\, M_x \cong M_y \,$ as (right) $A$-modules, there is an element
$\,\kappa \in Q \,$, unique up to a constant factor, such that $\, M_y = \kappa \,M_x \,$.
We can naturally extend $ \kappa $ to an {\it isomorphism of complexes}
$\, \Phi : \K_x \to \K_y \,$ making commutative the diagram
\[
\begin{diagram}[small, tight]
 M_x                 & \rTo^{r_x}           & \K_x       \\
 \dTo^{\kappa \cdot} &                      & \dTo_{\Phi}\\
 M_y                 & \rTo^{r_y}           & \K_y        \\
\end{diagram}
\]
To do this we need some extra notation. First, we denote by $\,
k(x)(y) \,$  (resp., $\, k(y)(x) \,$) the subspace of $\, Q \,$
spanned by elements of the form $\, f(x)\,g(y) \,$ (resp., $\,
g(y)\,f(x) \,$) with $\, f(x) \in k(x) \,$ and $\, g(y) \in
k(y)\,$. Next, we extend  (\ref{rx}) and (\ref{ry}) to these
subspaces. More precisely, we define the four linear maps:
\begin{equation}
\la{maps}
\begin{diagram}[small, tight]
        &                 & k(x)(y)  &                    &    \\
        & \ldTo^{\rrho_x} &   & \rdTo^{\lrho_y}   \\
k[x](y) &                 &   & & k(x)[y] \\
\end{diagram}
\qquad \quad
\begin{diagram}[small, tight]
        &                 & k(y)(x)  &                    &    \\
        & \ldTo^{\lrho_x} &   & \rdTo^{\rrho_y}   \\
k(y)[x] &                 &   & & k[y](x) \\
\end{diagram}
\end{equation}
where the accents indicate ``on which side'' the polynomial part
is taken.
For example, $\, \rrho_x: k(x)(y) \to k[x](y) \,$ is given by
$ f(x)\,g(y) \mapsto f(x)_{+}\,g(y)\,$.

Now, given a triple $\,(M_x,\,M_y, \,\kappa)\,$ as above, we define
$\, \phi: A \to A \,$ by
\begin{equation}
\la{phi}
\phi(a) := \rrho_y\,\lrho_x(\kappa \cdot a)\ , \quad a \in A\ .
\end{equation}
Note that (\ref{phi}) makes sense since $\, \kappa \in k(y)(x) \,$ and
$\, k(y)(x) \,$ is closed in $ Q $ under the right (and left) multiplication
by elements of $ A\,$.
\begin{lemma}
\la{L5}

$(1)\ $ $ \phi \,$ extends $\, \kappa \,$ through $\, r_x \,$,
i.~e. $\, \phi \circ r_x = r_y \circ \kappa \,$.

 $(2)\ $ $ \phi $ is invertible with $\, \phi^{-1}: A \to A \,$ given
by $ \phi^{-1}(a) = \rrho_x\,\lrho_y (\kappa^{-1} \cdot a) \,$.

$(3)\ $ We have $\, \phi(a) = a \,$ whenever $\, a \in k[x] $ or $ a \in k[y]\,$.
\end{lemma}
\begin{remark}
\la{R6}
In \cite{BW2} the map $\, \phi \,$ is denoted by $\Phi$ and defined by a different
formula (cf. \cite{BW2}, (5.4)). Lemma~\ref{L5}$(1)$ implies that the two definitions
in fact coincide.
\end{remark}
\begin{proof}
Denote by $\, k(x)_{-} \subset k(x) \,$ (resp.,  $\, k(y)_{-} \subset k(y) $)
the subspace of functions vanishing at infinity, so that
$\, k(x) = k[x] \oplus k(x)_{-}\,$ (resp., $\, k(y) = k[y] \oplus k(y)_{-}$).
Then we can extend our earlier notation writing, for example,
$\, k(y)_{-}(x)_{-} \,$ for the subspace of $ k(y)(x) $ spanned by all elements
$\, f(y)\,g(x) \,$ with $\, f(y) \in k(y)_{-}\,$ and $\, g(x) \in k(x)_{-}\,$.
With this notation, it is easy to see that $\, \kappa \in 1 +  k(y)_{-}(x)_{-} \,$
and $\, \kappa^{-1} \in 1 + k(x)_{-}(y)_{-} \,$ (cf. \cite{BW2}, Proposition~5.2(iii)).

$(1)\,$ Since $\, M_x \subset k(x)[y] \,$, $\, r_x(m) - m \in k(x)_{-}[y] =
k[y](x)_{-}\,$ for any $ m \in M_x\,$. Hence
$\, \kappa \cdot r_x(m) - \kappa \cdot m \in k(y)(x)_{-}\,$ and therefore
$\, \lrho_x(\kappa \cdot r_x(m)) = \lrho_x(\kappa \cdot m)\,$.
On the other hand, if $\, m \in M_x \,$ then $ \,\kappa \cdot m \in M_y \subset k(y)[x] $
and $ \lrho_x(\kappa \cdot m) = \kappa \cdot m $. Combining these together,
we get $\, \rrho_y\,\lrho_x(\kappa \cdot r_x(m)) =
\rrho_y(\kappa \cdot m) = r_y(\kappa \cdot m) \,$, which is equivalent to $ (1)\,$.

$(2)\,$ It follows trivially from (\ref{phi}) that
$\, \rrho_y\,\lrho_x(\phi(a) - \kappa \cdot a) = 0 \,$ for all $ a \in A $.
Since
$$
\Ker(\rrho_y\,\lrho_x) = (\lrho_x)^{-1}[\,\Ker\,\rrho_y\,] =
(\lrho_x)^{-1}[\,k(y)_{-}(x)\,] = k(y)(x)_{-} + k(y)_{-}(x)\ ,
$$
we have
$$
\phi(a) - \kappa \cdot a \, \in \, k(y)(x)_{-} + k(y)_{-}(x)  =
k[y](x)_{-} + k[x](y)_{-} + k(y)_{-}(x)_{-}\ .
$$
Multiplying this by $ \kappa^{-1} $ (and taking into account that
$\, \kappa^{-1} \in 1 + k(x)_{-}(y)_{-} $) yields
\begin{equation}
\la{inc2}
\kappa^{-1}\cdot \phi(a) - a \, \in \,
k[y](x)_{-} + k(x)(y)_{-} + k(y)_{-}(x)_{-} + k(x)_{-}(y)_{-}(x)_{-}\ .
\end{equation}
On the other hand, as $\, \phi(a) \in A\,$, we have
$\, \kappa^{-1} \cdot \phi(a) - a \in k(x)(y) \,$, so the expression
$\,\rrho_x\,\lrho_y (\kappa^{-1} \cdot \phi(a) - a) \,$ makes sense.
We claim that
\begin{equation}
\la{inc3}
\rrho_x\,\lrho_y (\kappa^{-1} \cdot \phi(a) - a) = 0 \quad \mbox{for all}\ a \in A\ .
\end{equation}
Indeed, fix $ a \in A $ and let
$\, \rrho_x\,\lrho_y (\kappa^{-1} \cdot \phi(a) - a) = b \,$ for some $\, b \in A \,$.
Then
$$
\kappa^{-1} \cdot \phi(a) - a - b \, \in \, \Ker(\rrho_x\,\lrho_y) =
k(x)_{-}(y) + k(x)(y)_{-} =  k[y](x)_{-} + k(x)(y)_{-}\ ,
$$
and hence, in view of \eqref{inc2}, $\, b \,$ can be written as a sum
$\, b = b_1 + b_2 \,$ with $\, b_1 \in k[y](x)_{-} \,$ and
$\, b_2 \in k(x)(y)_{-} + k(y)_{-}(x)_{-} + k(x)_{-}(y)_{-}(x)_{-} \,$.
Now, clearing denominators, we can find a polynomial $\, p \in k[x] \,$ such that
$ b_1 p \in A $ and $ b_2 p \in k(x)(y)_{-}\,$. Multiplying $\, b = b_1 + b_2 \,$ by
$\, p \,$ and applying $ \lrho_y $ to the resulting equation, we get
$\, b p = b_1 p \,$ which, in turn, implies the equality $\, b = b_1 \,$.
Since $\, b \in A \,$, $\, b_1 \in k[y](x)_{-} \,$ and $\, A \,\cap\, k[y](x)_{-} = \{0\}\,$,
we conclude  $\, b = 0 \,$, thus proving \eqref{inc3}.

It follows from \eqref{inc3} that $\, \rrho_x\,\lrho_y (\kappa^{-1} \cdot \phi(a)) = a \,$ for all
$ a \in A\,$. Defining now $\, \psi: A \to A \,$ by
$\, \psi(a) := \rrho_x\,\lrho_y (\kappa^{-1} \cdot a) \,$ we see that
$ \psi \circ \phi = \id_{A} \,$. On the other hand, reversing the roles
of $ \phi $ and $ \psi $ in the above argument gives
obviously $\,\phi \circ \psi = \id_{A} \,$. Thus, $ \phi $ is an isomorphism
of vector spaces, $ \psi $ being its inverse.

$(3)$ is immediate from the definition of $ \phi\,$. For example, if
$ a \in k[x] $  then $ \kappa\cdot a - a \in k(y)_{-}(x) $ and therefore
$ \phi(a) = \rrho_y\,\lrho_x(\kappa \cdot a) = \rrho_y\,\lrho_x(a) = a \,$,
as claimed.
\end{proof}

In view of Lemma~\ref{L5}, $\,\phi \,$ induces naturally the isomorphism of
quotient spaces $\, \bar{\phi}: V_x \to V_y \,$, and hence the isomorphism of complexes
$\, \Phi: \K_x \to \K_y \,$. We can use $ \Phi $ to identify $ \K_x $ and $ \K_y $ and
transfer the algebraic structure from one complex to another. More precisely, we set
$\, \K := \K_x \,$, i.~e. $\, K^0 := A \,$ and $\, K^1 := V_x \,$, and denote by
$\, m_1: K^0 \to K^1 \,$ the canonical projection. Next, we fix the ``cyclic'' vectors:
\begin{equation}
\la{cycl}
i:= 1 \in A = K^0 \ , \quad  \i := m_1(i) \in K^1 \ ,
\end{equation}
and define the endomorphisms
$\, X, \,Y \in \End_{k}(K^0) \,$ and $\, \X, \,\Y \in \End_{k}(K^1) \,$ by
\begin{eqnarray}
&& X(a) := \phi^{-1}(\phi(a)\cdot x)\ , \quad Y(a) := a \cdot y \ , \la{XY} \\*[1ex]
&& \X(\bar{a}) := \bar{\phi}^{-1}(\bar{\phi}(\bar{a})\cdot x)\ , \quad
\Y(\bar{a}) := \bar{a} \cdot y \ ,\la{qXY}
\end{eqnarray}
where ``$\,\cdot\,$'' stands for the usual multiplication in $ A $ and
$\, \bar{a} \in K^1 = V_x \,$ for the residue class of $\, a \in A \,$
mod $\,r_x(M_x)\,$. Clearly $\, \X \, m_1 = m_1 \, X \,$ and $\,\Y \, m_1 = m_1 \, Y \,$.
Moreover, we have the following crucial
\begin{proposition}
\la{P2}
The endomorphisms \eqref{XY} and \eqref{qXY} satisfy the equations
$$
X\,Y - Y \,X + \id_{K^0} = i\,j \ , \quad  \X\,\Y - \Y\,\X + \id_{K^1} = \i\, \j
$$
for some $\, j: K^0 \to k \,$ and $\, \j: K^1 \to k \,$ related by $\, j = \j\,m_1 \,$.
\end{proposition}
\begin{proof}
It suffices to show that
\begin{equation}
\la{const}
\left(X\,Y - Y\,X\right)a + a \,\in \,k \ , \quad \forall\, a \in A \ .
\end{equation}
Indeed, if \eqref{const} holds we may simply define $\, j(a) := \left(X\,Y - Y\,X\right)a + a \,$
satisfying the first equation of Proposition~\ref{P2}.
By Lemma~\ref{L5}$(1)\,$, it is then easy to see that
$\, j(a) = 0 \,$ on $\, \im(r_x) \,$, and since $ \im(r_x) = \Ker(m_1) $
the second equation follows from the first.

Now, to prove (\ref{const}) we start with equation \eqref{inc3} which is equivalent to
$$
\kappa^{-1} \cdot \phi(a) - a \, \in k(x)_{-}(y) + k(x)(y)_{-} = k[y](x)_{-} + k(x)(y)_{-}\ .
$$
Multiplying this by $ x $ on the right, we get
$$
\kappa^{-1} \cdot \phi(a) \cdot x - a \cdot x \, \in \, k[y] + k[y](x)_{-} + k(x)(y)_{-} =
k[y] + k(x)_{-}(y) + k(x)(y)_{-} \ ,
$$
whence the inclusion
$$
\rrho_x\,\lrho_y (\kappa^{-1} \cdot \phi(a) \cdot x - a \cdot x) =
\rrho_x\,\lrho_y (\kappa^{-1} \cdot \phi(a) \cdot x) - a \cdot x  \in  k[y]\ .
$$
By Lemma~\ref{L5}(2), this can be written as
$\,\phi^{-1}(\phi(a) \cdot x) - a \cdot x  \in  k[y] \,$, or equivalently
\begin{equation}
\la{rky}
X(a) - a \cdot x  \in  k[y] \quad \mbox{for all}\ a \in A \ .
\end{equation}
Now, using \eqref{rky}, we observe
\begin{equation}
\la{rk1y}
[X,\,Y]a + a = \left(X(a \cdot y) -
(a \cdot y) \cdot x\right) - \left(X(a) - a \cdot x\right) \cdot y
\, \in \, k[y]\ .
\end{equation}
On the other hand, define $\, X',\,Y' \in \End_{k}(A)\,$ by
$\, X' := \phi \, X\, \phi^{-1} $ and $\, Y' := \phi \, Y\, \phi^{-1} $,
so that
$$
X'(a) = a \cdot x \ ,\quad Y'(a) = \phi(\phi^{-1}(a) \cdot y)\ .
$$
Arguing as above, we can show then that $\, Y'(a) - a \cdot y \in k[x] \,$ for all
$\,a \in A\,$, which, in turn, yields the inclusion
$$
[X',\,Y']a + a = \left(Y'(a) - a \cdot y\right)\cdot x -
\left(Y'(a \cdot x) - (a \cdot x)\cdot y \right) \, \in \, k[x] \ .
$$
It follows now that $\, \phi\left([X,\,Y]a + a\right) = [X',\,Y']\phi(a) + \phi(a) \in k[x]\,$,
and therefore
\begin{equation}
\la{rk1x}
[X,\,Y]a + a \, \in \, \phi^{-1}(k[x])\ .
\end{equation}
By Lemma~\ref{L5}$(3)$, $\, \phi^{-1}(k[x]) = k[x] \,$, so comparing \eqref{rk1y} and
\eqref{rk1x}, we see that $\, [X,\,Y]a + a \, \in \, k[y]\,\cap \, k[x] = k\,$
as claimed in \eqref{const}.
\end{proof}

By Proposition~\ref{P2}, the complex $ \K $ together with linear
data $ (X, Y, i, j) $ and $ (\X, \Y, \i, \j) $  satisfies the conditions
\eqref{2.10} and \eqref{2.11}, and hence by Lemma~\ref{L6}, defines a homotopy
module over $A$. Using the lifting \eqref{lift1} we can now refine $ \K $
into an $ \A$-module as in Section~\ref{Recon}. The corresponding structure maps
\eqref{oper} satisfy the conditions \eqref{2.3}--\eqref{2.5} automatically.
The finiteness axiom \eqref{2.1} follows from \cite{BW2}, Proposition~5.2, and,
with our identification $ K^0 := A $ and the choice of cyclic vector \eqref{cycl},
$\, m_2^0(i,\,\mbox{--}\,): A \to K^0 \,$ becomes the identity map:
$$
m_2^0(i,\,x^k y^m) = \varrho^0(x^k y^m)i = Y^m X^k(i) = x^k y^m\ .
$$
This finishes our first construction of $\A$-envelopes.

\subsection{The second construction}
In this section $ A $ stands, as usual, for the Weyl algebra and $ A_0 $
for the commutative polynomial ring $ k[\bar{x},\bar{y}]  $ in variables $ \bar{x} $ and
$ \bar{y} $. We fix the lexicographic order on monomials of $ A $ and $ A_0 $ setting
$$
x^k y^l \prec x^{k'} y^{l'} \quad \mbox{and}\quad
\bar{x}^k \bar{y}^l \prec \bar{x}^{k'}\bar{y}^{l'} \quad \Longleftrightarrow \quad
l < l' \ \, \mbox{or}\ \, k< k'\  \mbox{if}\ l=l' \ .
$$
Given an element $\, a \in A \,$ (resp., $\, a \in A_0 $), we
write $ \sigma(a) \in A $ (resp., $\, \sigma(a) \in A_0 $) for the
initial (= greatest) term of $\,a\,$ with respect to this order,
and abbreviate ``{\it l.t.}'' for the lower terms $\,a -
\sigma(a)\,$.

Now, let $M$ be a rank $1$ torsion-free $A$-module. We fix an embedding $ M \into A $ and write
$\, \Sigma := \{(k,\,l) \in \mathbb N\times\mathbb N\,:\, \sigma(m) = x^ky^l \ \mbox{for some}\ m \in M \}\,$
for the set of exponents of $ M $ with respect to $ \prec \,$. Next, we set
$\, M_0 := \mbox{span}_k\{\bar{x}^k \bar{y}^l\,:\, (k,\,l) \in \Sigma\} \,$, which is obviously a
monomial ideal of $A_0$. Choosing an element $\, m = x^ky^l+l.t.\in M \,$,
one for each monomial $\, \bar{x}^k \bar{y}^l \in M_0 \,$, defines a linear isomorphism
$\, M_0 \stackrel{\sim}{\to} M \,$. We fix one of such isomorphisms and denote by
$\,r:\, M \to M_0 \,$ its inverse. The action of $ x $ and $ y $
on $M$ induces then two endomorphisms $X$ and $Y$ of $M_0$ satisfying
\begin{equation}
\la{in}
r(m.x)=X r(m)\quad \mbox{and}\quad r(m.y)=Y r(m)\quad\mbox{for all}\ m\in M\ ,
\end{equation}
and we have $\,X(a)=\bar x.a+ l.t.\,$
and $\,Y(a) =\bar y.a +l.t.\, $ for all $\, a \in M_0 \,$.

In general, $M_0\subseteq A_0$ has infinite codimension; however,
if we take the {\it minimal} principal ideal  $I\subseteq A_0$
containing $M_0$, then
\begin{equation}
\la{fini} \dim_k (I/M_0) < \infty\ .
\end{equation}
As an $A_0$-module, $I$ is free and generated by some monomial, which we denote $\,i\,$.
Next, we extend (somewhat arbitrarily) the endomorphisms $X$ and $Y$ from $M_0$ to $I$ by letting
$Xa:=\bar x.a$ and $Ya:=\bar y.a$ for all {\it monomials} $\,a=\bar x^k\bar y^l\in
I\setminus M_0\,$. The resulting maps still satisfy the properties
\begin{equation}
\la{tria} X(a) =\bar x.a+ l.t.\quad \mbox{and} \quad Y(a) =\bar y.a
+l.t. \quad\mbox{for any}\ a\in I\ ,
\end{equation}
and, as $\,xy-yx=1\,$, the following relation
\begin{equation}
\la{cr} [X,Y]+\id=0\quad\mbox{on}\quad M_0\subseteq I\,.
\end{equation}
Also, in view of \eqref{tria}, the elements $\,Y^lX^k(i)\,$ with $\,k,l \geq 0\,$ form a linear 
basis in $\,I\,$.

The above data satisfy the assumptions of the following
proposition which is crucial for our construction of $\A$-envelopes.
\begin{proposition}
\la{cf}
Let $ I $ be a rank $1$ free $ A_0$-module with generator
$\, i \,$ and induced order $ \prec $ (as defined above). Let $
\,M_0 \, $ be a subspace of $\, I \,$, stable under a pair of
linear endomorphisms $\, X, Y \in\End_k(I)\,$, satisfying
\eqref{fini}, \eqref{tria} and \eqref{cr}. Then there is another
pair $\, X, Y \in\End_k(I)\,$ that agrees with the given one
on $\,M_0\,$, satisfies the properties \eqref{tria}, \eqref{cr}
and, in addition,
\begin{equation}
\la{ro} 
\im ([X,Y]+\id) \subseteq k.i\,.
\end{equation}
\end{proposition}

Assuming (for the moment) that Proposition~\ref{cf} is true, we
complete our construction of an $\A$-envelope of $M$. To this end,
let $\, \K :=[\,0 \to K^0 \stackrel{m_1}{\longrightarrow} K^1 \to
0 \,]\,$ be a complex with $\,K^0 := I\,$, $\,K^1 := I/M_0 \,$ and
$m_1$ given by the canonical quotient map. Equip $ K_0 $ with
endomorphisms $X$ and $Y$ (granted by the above proposition) and
define a functional $\,j: K^0\to k\,$ by \eqref{ro} so that
$\,[X,Y]+\id_{K^0}=ij\,$. As $ m_1 $ is surjective, these maps 
induce linear maps $\X$, $\Y$ and $\j$ on $ K_1 $ satisfying
\eqref{2.8}, \eqref{2.10} and \eqref{2.11}. Thus, we obtain a
complex $\K$ of vector spaces with $\,(X,Y,i,j)\,$ and
$\,(\X,\Y,\i,\j)\,$, satisfying \eqref{2.11}, and a quasi-isomorphism $r: M
\into \K$ satisfying \eqref{in}. As shown in Section \ref{Recon},
these data determine an $\A$-envelope of $M$.

Now we turn to the proof of Proposition \ref{cf}. We will do this
in two steps: first, we will ``modify'' $X$ to achieve the inclusion
\begin{equation}
\la{r12} 
\im([X,Y]+\id) \subseteq k[\bar x].i\ ,
\end{equation}
then we will ``modify''  $Y$ to achieve \eqref{ro}. 

The first step is described by the following lemma. Note  
the condition $(b)$ of this lemma guarantees that the ``modification''
$\, (X,Y) \rightsquigarrow (X+X',Y) \,$ preserves the properties \eqref{tria}.
\begin{lemma}
\label{1}
Under the assumptions of Proposition \ref{cf}, there is
$\,X'\in\End_k(I)\,$ such that

$(a)\ $ $X'\equiv 0$ on $M_0\,$,

$(b)\ $ $X'(a) \prec \bar x.a$ for all $a\in I \,$,

$(c)\ $ $ \im([X+X',Y]+\id)\subseteq k[\bar x].i\,$.
\end{lemma}
\begin{proof} 
Since $\,M_0 \,$ is a subspace of finite codimension in $I$, 
invariant under the action of $ Y $, there is a filtration on $I$:
$\,M_0 \subset M_1\subset \ldots \subset M_n=I\,$, such that
$\,YM_j\subseteq M_j\,$ and $\,\dim_k(M_j/M_{j-1})=1 \,$. Choose
$\,v_1, v_2,  \ldots, v_n \in I\,$, so that $\,M_j=M_{j-1}\oplus k.v_j\,$.
If $\,\alpha_j\in k\,$ are the eigenvalues of the maps induced by 
$Y$ on the quotients  $\,M_j/M_{j-1}\,$, we have
\begin{equation}
\label{tri}
m_{j-1}:=Y(v_j) - \alpha_jv_j\in M_{j-1}\ ,\quad j=1,2,\ldots, n \ .
\end{equation}
Now, we will construct $\, X'\in\mathrm{End}_k(I)\, $ by setting
$\,X'\equiv 0\,$ on $\,M_0\,$ and defining $\,X'(v_j)\,$
successively for $\,j=1,2,\ldots, n\,$. At each step, we will
verify that
\begin{equation}
\la{r11}
([X+X',Y]+\id)v_j \in k[\bar x].i\ .
\end{equation}
Clearly, the last condition of Lemma~\ref{1} follows from \eqref{r11}.

Suppose that we have already defined $X'$ on $M_{j-1}$, and it
satisfies the condition $(b)$ of the lemma. (This is obviously the
case for $j=1$.) A trivial calculation then shows
\begin{equation}
\la{com} ([X+X',Y]+\id)v_j = ([X,Y]+\id)v_j +
X'(m_{j-1}+\alpha_jv_j)-YX'(v_j)\ ,
\end{equation}
where $\,m_{j-1}\in M_{j-1} \,$ is given by \eqref{tri} and $\,X'(v_j)\,$ 
has yet to be defined.

By our induction assumption, the expression $\,([X,Y]+\id)v_j +
X'(m_{j-1})\,$ is already defined, and we denote it by $\,u\,$.
The right-hand side of \eqref{com} then becomes
$\,u-(Y-\alpha_j)X'(v_j)\,$. Now, to satisfy \eqref{r11}, it suffices 
to find $\,a\in I\,$ such that $\,u-(Y-\alpha_j)a \in k[\bar
x].i\,$. To this end, using \eqref{tria}, one can show easily that
every $\,u \in I\,$ can be written as
\begin{equation}
\la{re} u = (Y-\alpha_j)a+b \quad\mbox{for some}\quad a\in I
\quad\mbox{and}\quad b\in k[\bar x].i\ .
\end{equation}
Thus, if $\, u = ([X,Y]+\id)v_j + X'(m_{j-1}) \,$, we take $\,a
\in I \,$ as in \eqref{re} and let $\,X'(v_j):= a \,$. Then
\eqref{r11} follows from \eqref{com}.

Finally, we check that the condition $ (b) $ holds on each 
filtration component $\,M_{j}\,$. By induction assumption, 
we have $\,X'(m_{j-1})\prec \bar x.m_{j-1}= (\bar
x\bar y).v_j+l.t.\,$; whence $\,X'(m_{j-1})\prec (\bar x\bar
y).v_j\,$. Furthermore, it follows from \eqref{tria}  that $\,
([X,Y]+\id)v_j \prec (\bar x\bar y).v_j\,$. Thus we have $\, u
\prec (\bar x\bar y).v_j\,$. On the other hand,  \eqref{re}
implies that $\,u=\bar y.a+l.t.\,$. Combining these last two
facts, we see that $\,X'(v_j)= a \prec \bar x.v_j\,$. So, if
the condition $(b)$ holds on $ M_{j-1} $, then it holds also on
$\,M_j=M_{j-1}\oplus k.v_j\,$. This completes the induction and
the proof of our lemma.
\end{proof}

Now, assuming \eqref{r12}, we will ``modify'' $\,Y
\rightsquigarrow Y + Y' \,$, so that the new endomorphisms $\,
(X,Y)\,$ satisfy \eqref{ro}. Again, the condition $(b)$ of
Lemma~\ref{2} below guarantees that \eqref{tria} remains true
after such a ``modification.''
\begin{lemma}
\label{2} 
In addition to the assumptions of
Proposition \ref{cf}, suppose that $\,X,Y\,$ satisfy \eqref{r12}.
Then there is $\,Y'\in\End_k(I)\,$ such that

$(a)\ $ $Y'\equiv 0$ on $M_0$,

$(b)\ $ $\im(Y') \subseteq k[\bar x].i$,

$(c)\ $ $ \im([X,Y+Y']+\id)\subseteq k.i$.
\end{lemma}

\begin{proof}
We will argue as in the proof of Lemma \ref{1}. We start by fixing
a  filtration $\,M_0\subset M_1\subset\dots \subset M_n=I\,$
on $\,I\,$, stable under the action of $ X $ and such that
$\,\dim_k(M_j/M_{j-1})=1\,$ for each $\,j = 1,2, \ldots, n\,$.
Then we choose $\,w_1, w_2, \ldots, w_n\in I \,$ so that
$\,M_j=M_{j-1}\oplus k.w_j\,$ and define the elements
\begin{equation}
\label{tri1}
m_{j-1}:=X(w_j) -\beta_j w_j\in M_{j-1}\ ,\quad j
=1,2,\ldots, n\ .
\end{equation}
where $\,\beta_j \in k \,$ are the eigenvalues of the maps induced by 
$X$ on the quotients $\,M_j/M_{j-1}\,$.

Next, setting $\,Y'\equiv 0\,$ on $\,M_0\,$, we will define
$\,Y'(w_j)\,$ successively for $\,j=1,2,\dots, n\,$, so that
$\,Y'(w_j)\in k[\bar x].i\,$ and
\begin{equation}
\la{r13} ([X,Y+Y']+\id)w_j \in k.i\ .
\end{equation}

Suppose that $\,Y'\,$ is already defined on $\,M_{j-1}\,$ and
$\,Y'(m) \in k[\bar x].i \,$ for all $\, m \in M_{j-1}\,$. (This
is obviously true for $\,j=1\,$.) Then, using \eqref{tri1}, we can write
\begin{equation}\la{com1}
([X,Y+Y']+\id)w_j=([X,Y]+\id)w_j - Y'(m_{j-1})+(X-\beta_j)Y'(w_j) \ .
\end{equation}
Note that in view of \eqref{r12}, $\,([X,Y]+\id)w_j \in k[\bar x].i\,$, and $\,
Y'(m_{j-1}) \in k[\bar x].i\,$ by our induction assumption. Hence 
$\, u:= ([X,Y]+\id)w_j - Y'(m_{j-1}) \in k[\bar x].i\,$. The right-hand side 
of \eqref{com1} then becomes
$\,u+(X-\beta_j)Y'(w_j)\,$. So, given $\,u\in k[\bar x].i\,$, it
suffices to show that there exists $\,a\in k[\bar x].i\,$ such that $\,
u+(X-\beta_j)a \in k.i\,$. But this follows easily from \eqref{tria}. 
Taking such an element $ \,a\,$ for $\, u= ([X,Y]+\id)w_j - Y'(m_{j-1}) \,$ 
and letting $\,Y'(w_j) := a\,$, we get \eqref{r13} as a consequence of \eqref{com1}. This
finishes the induction and the proof of Lemma~\ref{2}, as well as the
proof of Proposition \ref{cf}.
\end{proof}

\section{The Calogero-Moser Correspondence}
\label{CMC}
Let $ \MM $ be the set of {\it strict} isomorphism classes of $ \A$-modules satisfying the
axioms \eqref{2.1}--\eqref{2.5}. In this section we establish two natural bijections between
$ \MM $ and $\,(a)\,$ the set $ \RR $ of isomorphism classes of (right) ideals in $ A \,$,
$\,(b)\,$ the (disjoint) union $ \CC $ of Calogero-Moser spaces $ \CC_n$, $\, n\geq 0\,$.
Combining these bijections we then recover the one-one correspondence $\,\RR  \leftrightarrow \CC\,$
constructed in \cite{BW1,BW2}.
\begin{theorem}
\la{T6}
There are four maps
\begin{equation}
\la{D7}
\begin{diagram}[small, tight]
\RR &  \pile{\rTo^{\theta_1}\\ \lTo_{\omega_1}} & \MM & \pile{\rTo^{\theta_2}\\ \lTo_{\omega_2}} & \CC
\ , \\
\end{diagram}
\end{equation}
such that $ (\theta_1, \omega_1) $ and $ (\theta_2, \omega_2) $ are pairs of mutually inverse
bijections, $\, \omega_1 \circ \omega_2  \,$ is the map $ \omega $ defined in \cite{BW1}
and $ \theta_2 \circ \theta_1  $ is the inverse of $ \omega $ constructed in \cite{BW2}.
\end{theorem}
\begin{proof}
All the maps have already been defined (implicitly) in the previous sections.

First, $\, \theta_1 $ is given by the constructions of Section~\ref{Ex} which assigns
to an ideal $ M $ its $ \A$-envelope $ M \stackrel{r}{\to} \K $. Since passing from $ M $
to an isomorphic module results only in changing the quasi-isomorphism $\, r \,$ (not $ \K $),
this indeed gives a well-defined map from $ \RR $ to $ \MM\,$.

Second, $\, \omega_1 \,$ is defined by taking cohomology of $\A$-modules:
$\, [\K] \mapsto [H^0(\K)] \,$. By Lemma~\ref{L3.1}, this makes sense since
$ H^0(\K) $ is a f.~g. rank $1$ torsion-free module over $ A $ and hence its isomorphism
class is indeed in $ \RR$. With this definition the equation
$\, \omega_1 \circ \theta_1 = \id_{\RR} \,$ is obvious while
$\, \theta_1 \circ \omega_1 = \id_{\MM} \,$ follows immediately from Theorem~\ref{T3}.
The maps $ \theta_1 $ and $ \omega_1 $ are thus mutually inverse bijections.

Third, in Section~\ref{CMD} we have shown how to obtain the Calogero-Moser data from an
$\A$-module $ \K $ satisfying \eqref{2.1}--\eqref{2.5}. Specifically,
we associate to $ \K $ the pair of endomorphisms $ (\X, \Y) $ arising from the
action of $x$ and $y$ on $ K^1\,$ together with a cyclic vector $ \i \in K^1 $
and a covector $\, \j: K^1 \to k\,$ (see \eqref{2.9}, \eqref{2.10}).
By Lemma~\ref{L3}, these satisfy the relation
$\,[\X, \Y] + \id_{K^1} = \i\,\j \,$ and hence represent a point in $ \CC\,$.
A strict isomorphism of $\A$-modules commutes with the action
of $ A $ and hence transforms the quadruple $ (\X, \Y, \i, \j) $ into an equivalent one.
Thus, strictly isomorphic $\A$-modules yield one and the same point in $ \CC $, and
we get a well-defined map $\, \theta_2: \MM \to \CC  \,$.

Fourth, in Section~\ref{Recon}, starting with Calogero-Moser data $ (\X, \Y, \i,\j) $
we construct an $ \A$-module $ \K $ that satisfy \eqref{2.1}--\eqref{2.5}.
If we replace now $ (\X, \Y, \i,\j) $ by equivalent data,
then the functional \eqref{fun} remains the same, and hence so do the ideal \eqref{id}
and the $ R$-module $ K^0 $. On the other hand, the differential $ m_1 $ gets changed to
$\, g\,m_1 \,$. As a result, we obtain an $ \A$-module $ \KK $ strictly isomorphic
to $ \K \,$, the isomorphism $ \K \to \KK $ being given by $\, (\id_{K^0}, g) \,$.
Thus, the  construction of Section~\ref{Recon} yields a well-defined map
$ \omega_2: \CC \to \MM\,$.

Finally, it remains to see that $ \omega_2 $ and $ \theta_2 $ are mutually inverse bijections.
First of all, the composition $\, \theta_2 \circ \omega_2 \,$ being the identity on $ \CC $
is an immediate consequence of definitions. On the other hand,
$\,\omega_2 \circ \theta_2 = \id_{\MM} \,$ can be deduced from Theorem~\ref{T3}.
Indeed, if we start with an $ \A$-module $ \K $ and let $ \KK $ represent the
class $\, \omega_2 \,\theta_2[\K] \,$ then $ \K $ and $ \KK $ have equivalent
finite-dimensional data $ (\X, \Y, \i,\j) \,$, and hence, in view of \eqref{numb},
determine the same functionals \eqref{invar}. It follows now from (the implication
$\,(d) \Rightarrow (a)\,$ of) Theorem~\ref{T3} that $ [\KK] = [\K] $ in $ \MM\,$.
\end{proof}

Combined with Corollary~\ref{C1}, Theorem~\ref{T6} leads to an explicit description of
ideals of $ A $ in terms of Calogero-Moser matrices.
\begin{corollary}
\la{C3}
The map $\, \omega: \CC \to \RR \,$ assigns to a point of $\, \CC \,$ represented by
$\, (\X, \Y, \i, \j) \,$ the isomorphism class of the fractional ideals
\begin{eqnarray}
&& M_x := \det(\X - x)\,A + \kappa^{-1} \,\det(\Y - y)\,A \ , \nonumber \\*[1ex]
&& M_y := \det(\Y - y)\,A + \kappa \,\det(\X - x)\,A\ , \nonumber
\end{eqnarray}
where $\,\kappa \,$ is given by \eqref{kappa}.
\end{corollary}

\section{DG-Envelopes}
\la{HPT}
The most interesting feature of the Calogero-Moser correspondence is its equivariance with respect
to the Weyl algebra automorphism group $ G $. Both in \cite{BW1} and \cite{BW2}, this result was proved
in a rather sophisticated and roundabout way. The main problem was that the bijection
$\, \RR  \to \CC\,$ was defined in \cite{BW1, BW2} indirectly, by passing through a third
space\footnote{namely, the adelic Grassmannian $ \mbox{Gr}^{\mbox{\scriptsize \rm ad}} $ in \cite{BW1} and the
moduli spaces of rank $1$ torsion-free sheaves over a noncommutative $ \mathbb{P}^2 $ in \cite{BW2}.},
and the action of $ G $ on that space was difficult to describe.
Unfortunately, our Theorem~\ref{T6} has the same disadvantage: the axiomatics of $\A$-envelopes
(specifically, the axioms \eqref{2.3}--\eqref{2.5}) are not invariant under the action of $ G $,
and thus it is not obvious how $ G $ acts on $ \MM $.
In this section we resolve this problem in a simple and natural way. The key idea inspired
by \cite{Q2} is to replace $A$ by its DG-algebra extension $ \AB $ and work with
DG-modules over $ \AB $ instead of $\A$-modules over $A$.

\subsection{Axioms}
Let $ \AB $ denote the graded associative algebra $\, I \oplus R \,$
having two nonzero components:
the free algebra $ R = k\langle x,y\rangle $ in degree $0$ and the (two-sided) ideal
$ I := RwR $ in degree $ -1\,$. The differential on $ \AB $ is defined by the natural inclusion
$\, d : I \hookrightarrow R \,$ (so that $\,dw = xy-yx-1 \in R \,$ and $\,da \equiv 0\,$ for all
$\,a \in R\,$). Regarding $ R $ and $ A $ as DG-algebras (with single component in degree $0$)
we have two DG-algebra maps: the canonical inclusion $ \iota: R \to \AB $ and projection
$ \eta: \AB \to A \,$. The latter is a quasi-isomorphism of complexes which we can interpret,
following \cite{Q2}, as a ``length one resolution'' of $ A $ in the category of associative DG-algebras.
Now, if $ \DGMod(\AB) $ is the category of (right unital) DG-modules over $ \AB $ we have two
restriction functors $ \iota_{*}: \DGMod(\AB) \to \Com(R) $ and $ \eta_{*}:
\Com(A) \to \DGMod(\AB) \,$, each being an exact embedding. We may (and often will) identify
the domains of these functors with their images, thus thinking of $ \Com(A) $ as a full
subcategory of $ \DGMod(\AB) $ and $ \DGMod(\AB) $ as a subcategory of $ \Com(R)\,$.
Note that, under the first identification, a DG-module $\, \LB \in \DGMod(\AB) \,$ belongs to
$ \Com(A) $ if and only if the element $\, w \in \AB \,$ acts trivially on $ \LB \,$.

Let $ M $ be, as usual, a rank $1$ finitely generated torsion-free module over $A\,$.
\begin{definition}
\la{D2}
A {\it DG-envelope}\, of $ M $ is a quasi-isomorphism $\, q: M \to \LB \,$ in $ \DGMod(\AB) $,
where  $\, \LB = L^0 \oplus L^1 \, $ is a DG-module with two nonzero components
(in degrees $ 0 $ and $1$) satisfying the conditions:
\begin{itemize}
\item {\it Finiteness:}
\begin{equation}
\la{7.1}
\dim_{k}\,L^1  < \infty \ .
\end{equation}
\item {\it Existence of a cyclic vector:}
\begin{equation}
\la{7.2}
L^0 \ \mbox{is a cyclic $R$-module with cyclic vector $ i $}.
\end{equation}
\item {\it ``Rank one'' condition:}
\begin{equation}
\la{7.3}
\LB.w \subseteq  k.i\ ,
\end{equation}
where $ \LB.w $ denotes the action of $ w $ on $ \LB\,$.
\end{itemize}
\end{definition}
\noindent
The following properties are almost immediate from the above definition.

1. The differential on $ \LB $ is given by a surjective $R$-linear map:
$ d_{\LB}:  L^0 \to L^1 $. Indeed, $ d_{\LB}$ being surjective follows from \eqref{7.1}
(and the fact that $ A $ has no nontrivial finite-dimensional modules);
$ d_{\LB} $ commuting with action of $ R $  is a consequence of the Leibniz rule.
Together with \eqref{7.2} these two properties imply that $ L^1 $ is a cyclic $ R$-module generated by
$ \i := d_{\LB}(i) \in L^1 \,$. In particular, we have $ \i \not= 0 $ (unless $ L^1 = 0\,$).

2. Being of negative degree in $ \AB $, the element $\, w\,$ acts trivially on $ L^0 $.
Hence, \eqref{7.3} is equivalent, in effect, to $\, L^1.w \subseteq k.i \,$.
Defining now $\, \j: L^1 \to k \,$ by
\begin{equation}
\la{7.3e}
v.w = \j(v)i\ , \quad v \in L^1\ ,
\end{equation}
and setting $\, j := \j\,d_{\LB}: L^0 \to k \,$, we have our usual relations (cf. \eqref{2.11})
\begin{equation}
\la{7.4}
X\,Y - Y \,X + \id_{L^0} = i\,j \ , \quad  \X\,\Y - \Y\,\X + \id_{L^1} = \i\, \j\ ,
\end{equation}
where $ (X, \X) \in \End_{k}(\LB) $ and $ (Y, \Y) \in \End_{k}(\LB) $ come from the action of
$ x $ and $ y $ on the corresponding components of $ \LB $. Note in this case the equations
\eqref{7.4} arise from the Leibniz rule: the first one can be obtained by
differentiating the obvious identity $\, u.w = 0 \,$, $\,\forall u \in L^0\,$:
$$
0 =  (du).w + u.dw = \j(du)\,i + u.(xy-yx-1) = i\,j(u) + (YX-XY-\id_{L^0})u \ ,
$$
while the second by differentiating \eqref{7.3e}:
$$
\i\,\j(v)  = d(v.w) = - v.dw = - v.(xy-yx-1) = (\X\Y - \Y\X + \id_{L^1})v   \ .
$$

3. Define $\, \varepsilon: R \to k \,$  by $\, a \mapsto \j(\i.a)\,$
(cf. \eqref{fun}) so that $\,\i.(aw) = \varepsilon(a) \,i \,$ for all $ a \in R$.
Differentiating then the identity $\, i.aw = 0 \,$ yields
$$
i.\left[a(xy-yx-1) + \varepsilon(a)\right] = 0 \ , \quad \forall a \in R\ .
$$
Hence the multiplication-action map $\, R \to L^0 \,$, $ a \mapsto i.a \,$,
factors through the canonical projection $ R \onto R/J \,$, where $\,
J := \sum_{a \in R} \left(\,a(xy-yx-1) + \varepsilon(a)\,\right)R \,$ (cf. \eqref{id}).
We claim that the resulting map
\begin{equation}
\la{isom}
 \phi:\, R/J \to L^0 \ \mbox{is an isomorphism of $ R$-modules.}
\end{equation}
Indeed, in view of \eqref{7.2}, $ \phi $ is surjective.
On the other hand, by the Snake Lemma, the kernel of $ \phi $ coincides with
the kernel of the natural map: $\, \Ker(d_{\LB}\circ \phi) \to \Ker(d_{\LB})\,$.
Since both $ \Ker(d_{\LB}\circ \phi) $ and $ \Ker(d_{\LB}) $ are rank $1$ torsion-free
$A$-modules, the last map is injective, and hence so is $ \phi $.

\begin{theorem}
\la{T7}
Every rank $1$ torsion-free $A$-module has a DG-envelope in $ \DGMod(\AB) $ satisfying
the axioms of Definition~\ref{D2}.
\end{theorem}
\begin{proof}
We define a DG-module $ \LB $ together with quasi-isomorphism
$\, q: M \to \LB \,$ by reinterpreting the construction of Section~\ref{Ex}.
We let $ \K_x $ (see \eqref{Kxy}) be the underlying complex for $ \LB \,$,
but instead of giving $ \K_x $ the structure of a homotopy module over $A$, we make it
a DG-module over $\AB\,$. Specifically, using the notation \eqref{cycl}--\eqref{qXY},
we define the action of $ \AB $ on $ \K_x $ by
\begin{equation}
\la{dgact}
(u,v).x := (X(u), \X(v))\ ,\ (u,v).y := (Y(u), \Y(v))\ ,\ (u,v).w := (\j(v)i,\, 0) \ ,
\end{equation}
where $\, (u,v) \in \K_x \,$. The compatibility of \eqref{dgact} with Leibniz's rule amounts to
verifying the relations \eqref{7.4}, and this has already been done in Proposition~\ref{P2}.
With this definition of $ \LB $ the conditions \eqref{7.1}--\eqref{7.3} are immediate.
Setting $ q := r_x $ gives thus a required DG-envelope of $M$.
\end{proof}

\subsection{DG-envelopes vs. $\A$-envelopes} It is clear from the above discussion that the axiomatics
of DG-envelopes is closely related to that of $\A$-envelopes. To make this relation precise we will
view $ A $ and $ \AB $ as $ \A$-algebras\footnote{For basic definitions concerning $\A$-algebras
we refer the reader to Section~3 of \cite{K1}.}
and use the following simple observation to relate the corresponding categories of modules.
\begin{lemma}
\la{L10}
Let $ A $ be an associative $k$-algebra, and let $ \eta: R \onto A $ be an algebra
extension with augmentation ideal $ I := \Ker(\eta) \,$. Form the DG-algebra
$\, \AB := I \oplus R \,$ with differential $ d $ given by the natural inclusion
$\, I \into R \,$. Then, choosing a linear section $\, \varrho: A \to R \,$ of
$\, \eta \,$ is equivalent to defining a quasi-isomorphism of $\A$-algebras
$\, \vr : A \to \AB \,$ with $\, \eta \circ \vr = \id_{A}\,$.
\end{lemma}
\begin{proof}
By degree considerations, any $\A$-algebra morphism $ \vr: A \to \AB $ may have only two
nonzero components, $\, \vr_1: A \to R \,$ and $ \vr_2: A \otimes A \to I \,$, which satisfy
the single relation (cf. \cite{K1}, Sect.~3.4):
\begin{equation}
\la{vr}
\vr_1(ab) = \vr_1(a)\,\vr_1(b) + d \vr_2(a,b)\ , \quad \forall a,b \in R\ .
\end{equation}
Now, if  $\, \varrho: A \to R \,$ is a linear map such that $ \eta \circ \varrho = \id_{A} $,
its curvature $ \omega $ (see \eqref{curv}) takes values in $\, I \subseteq R \,$.
Hence, setting $ \vr_1 := \varrho $ and $ \vr_2 := d^{-1} \omega $ makes sense and obviously
satisfies  \eqref{vr}.
Conversely, if $ \vr: A \to \AB $ is a morphism of $ \A$-algebras with
(left) inverse $ \eta $ then, by definition, $ \vr_1: A \to R $ is a linear section of $ \eta $ and
$ d \vr_2 $ is its curvature.
\end{proof}

\begin{remark}
The above lemma is essentially due to Quillen: in \cite{Q2}, Section~5.1, he formulates
this result in a slightly greater generality using the language of DG-coalgebras.
\end{remark}

Given a morphism of $\A$-algebras $\, f: A \to B \,$, there is a natural restriction
functor $\, f_{*}: \Modi(B) \to \Modi(A) \,$ from the category of $\A$-modules over $ B $
to the category of $\A$-modules over $ A$. Specifically, $ f_{*} $ assigns to an $\A$-module
$\, \LB \in \Modi(B) $ with structure maps $\, m_n: \, \LB \otimes B^{\otimes(n-1)} \to
\LB \,$, $\,n\geq 1\,$, an $\A$-module $ f_{*}\LB \in \Modi(A) $ with the same underlying
vector space as $ \LB $ and structure maps $ m_{n}^{A}: \, \LB \otimes A^{\otimes(n-1)} \to \LB $
given by (see \cite{K1}, Section~6.2)
\begin{equation}
\la{mm}
m_{n}^{A} : = \sum\,(-1)^s m_{r+1}(\id \otimes f_{i_{1}} \otimes \ldots \otimes f_{i_{r}})\ ,
\quad n \geq 1\ .
\end{equation}
If $ \varphi: \LB \to \M $ is a morphism of $ \A $-modules over $ B $ the map
$ \, f_{*}\varphi: \, f_{*}\LB \to f_{*}\M \in \Modi(A) \,$ is defined by
\begin{equation}
\la{ff}
(f_{*}\varphi)_n : = \sum\,(-1)^s \varphi_{r+1}(\id \otimes
f_{i_{1}} \otimes \ldots \otimes f_{i_{r}})\ , \quad n \geq 1\ .
\end{equation}
The sums in \eqref{mm} and \eqref{ff} run over all $\, r: 1 \leq r \leq n-1 \,$
and all integer decompositions $\, n-1 = i_1 + \ldots + i_r \,$ with $\, i_k \geq 1 \,$, and
$\, s := \sum_{k=1}^{r-1} k(i_{r-k} - 1 )\, $.

If we apply this construction in the situation of Lemma~\ref{L10} (for $\, f = \vr $)
and take into account that $ \DGMod(\AB) $ is naturally a subcategory of $ \Modi(\AB) \,$, we
get a faithful functor $ \vr_{*}:\, \DGMod(\AB) \to \Modi(A) $ transforming DG-modules over
$ \AB $ to $ \A$-modules over $A$. Using this functor, we can state the
precise relation between $\A$- and DG-envelopes.
\begin{theorem}
\la{T12}
Let $ \varrho: A \to R $ be given by
\begin{equation}
\la{sect}
\varrho(x^k y^m) = x^k y^m\ ,\quad \forall\, k,m \geq 0\ .
\end{equation}
The corresponding restriction functor $ \vr_{*}: \DGMod(\AB) \to \Modi(A) $
gives an equivalence between the full subcategory of $ \DGMod(\AB) $
consisting of DG-modules that satisfy Definition~\ref{D2}
and the subcategory of $\A$-modules with {\rm strict} morphisms
satisfying  Definition~\ref{D1}. Under this equivalence a DG-envelope
$ q: M \to \LB $ of $ M $ transforms to its $\A$-envelope
$\, \vr_{*}q: M \to \vr_{*}\LB\,$.
\end{theorem}
\begin{proof}
First, using \eqref{mm} we compute the structure maps on $ \K := \vr_{*}\LB $
when $ \LB $ has only two nonzero components (in degrees $0$ and $1$):
\begin{eqnarray}
&& m_1(\xx) = d_{\LB}(\xx)\ ,\quad \xx \in \K \ , \la{m1} \\
&& m_2(\xx, a) = \xx.\vr_1(a)\ ,\quad \xx \in \K \ , \ a \in A \ ,\la{m2} \\
&& m_3(\xx, a, b) = \xx. \vr_2(a,b)\ ,\quad \xx \in \K \ , \ a, b \in A \ ,\la{m3}
\end{eqnarray}
and $ m_n  \equiv 0 $ for all $ n \geq 4 \,$.
Next, we observe that if $ \varrho $ is given by \eqref{sect} then
$ \omega(x, x^k y^m) = \omega(x^k y^m,y) = 0 \,$ for all $\,k, m\geq 0 \,$,
and $ \omega(y,x) = xy-yx-1 $. As $ \omega = d \vr_2 $ (see Lemma~\ref{L10})
and $ d: I \to R $ is injective, this implies
\begin{equation}
\la{rho2}
\vr_2(x, a) = \vr_2(a,y) = 0 \ , \quad \forall\,a \in A\ , \ \mbox{and}\ \vr_2(y,x) = w\ .
\end{equation}
Looking now at \eqref{m3} we see that \eqref{2.3} and \eqref{2.4} hold
automatically for $ \K $, while \eqref{2.5} follows from  \eqref{7.3}.

To verify \eqref{2.2}, we factor $\, m_2(i,\, \mbox{--}\,):\, A \to K^0,\,
a \mapsto i.\varrho(a) \,$, (see \eqref{m2}) as
$$
A \stackrel{\varrho}{\longrightarrow} R \onto R/J
\stackrel{\phi}{\longrightarrow}  L^0 = K^0\ ,
$$
where $ \phi $ is the quotient of the multiplication-action map $ R \to L^0 \,$,
$\,a \mapsto i.a \,$. If $ \varrho $ is defined by \eqref{sect} then the composition
of the first two arrows is obviously an isomorphism. On the other hand, if $ \LB $
satisfies Definition~\ref{D2} (specifically, $ H^0(\LB) $ is a rank $1$ torsion-free
$A$-module) then $ \phi $ is also an isomorphism (see \eqref{isom}), and hence so is $ m_2(i,\, \mbox{--}\,)$.
Finally, the axiom \eqref{2.1} is equivalent to \eqref{7.1}.

Thus, we have shown that $ \vr_* $ takes a DG-envelope of any $ M $
(in the sense of Definition~\ref{D2}) to its $\A$-envelope (in the sense of Definition~\ref{D1}).
Moreover, it is clear that {\it every} $\A$-envelope $\K$ satisfying \eqref{2.1}--\eqref{2.5} is
strictly isomorphic to one of the form $ \vr_* \LB $.
Indeed, given such a $\K$ we can define $ \LB  $ by identifying it with $ \K $
as (a complex of) vector spaces, and imposing a DG-structure on $\K$ as in
Theorem~\ref{T7} (see \eqref{dgact}).

It remains to see that $\,\vr_* \,$ transforms DG-morphisms to {\it strict} \,
$\A$-morphisms, and that every strict morphism between $\A$-envelopes
can be obtained in this way.
The first statement follows directly from formula \eqref{ff} which simplifies
in our case to $\,(\vr_* \varphi)_1 = \varphi_1 \,$,
$\,(\vr_{*}\varphi)_2 = \varphi_2(\id \otimes \varrho ) \,$ and
$ \,(\vr_{*}\varphi)_n = 0 $ for all $ n \geq 3\,$. (So if $ \varphi \in \DGMod(\AB) $
then $ \varphi_2 = 0 $ and hence $ \,(\vr_{*}\varphi)_n = 0 $ for all $ n \geq 2\,$.)
Conversely, if $ f: \K \to \KK $ is a strict morphism of $\A$-modules in
$ \Modi(A) $ then $ f = f_1 $ commutes both with differentials and the action of $x$ and $y$.
Moreover, we have $\, m_3(f^1(v),a,b) = f^0(m_3(v,a,b)) \,$ (see \eqref{1.11e})
which is equivalent, in view of \eqref{m3}, to $\, f^1(v). \vr_2(a,b) = f^0(v.\vr_2(a,b)) \,$.
Letting $ a=y $ and $ b = x $ and taking into account \eqref{rho2} we see that $ f $ also
commutes with $ w $, and hence with any element of $\AB$ (as $\AB $ is generated by $x,y,w $).
This shows that $ f $ is indeed a morphism of DG-modules over $ \AB $, thus finishing the proof of
the Theorem.
\end{proof}

\begin{corollary}
\la{U2}
A DG-envelope of a given module $M$ is determined uniquely, up to unique isomorphism in $\DGMod(\AB)$.
\end{corollary}
\begin{proof}
Combine Theorem~\ref{T12} with Corollary~\ref{U1}.
\end{proof}
\begin{corollary}
\la{bij}
There is a natural bijection $ \varrho_* $ between the set $ \tilde{\MM} $ of isomorphism classes of DG-modules satisfying
Definition~\ref{D2} and the set  $ \MM $ of strict isomorphism classes of $\A$-modules satisfying
Definition~\ref{D1}.
\end{corollary}
\begin{proof}
The bijection $ \varrho_* $ is induced by the equivalence of Theorem~\ref{T12}.
\end{proof}

\subsection{G-equivariance}
Let $ G $ be the group of $k$-linear automorphisms of the free algebra $R$ preserving the
commutator $\,xy-yx \in R\,$. Every $ \sigma \in G $ extends uniquely to an automorphism
$ \tilde{\sigma} $ of the graded algebra $ \AB \,$ in such a way that $\,\tilde{\sigma}\,d_{\AB} =
d_{\AB}\,\tilde{\sigma}\,$. Specifically, $\, \tilde{\sigma}: \AB \to \AB \,$ is defined on generators
by $\, \tilde{\sigma}(x) = \sigma(x) \,$, $\,\tilde{\sigma}(y) = \sigma(y) $ and $\tilde{\sigma}(w) = w \,$.
Thus, we have an embedding $\, G \into \DGAut_k(\AB)\,$, where $ \DGAut_k(\AB) $ is the group of
all DG-algebra automorphisms of $\AB$.

For any $ \sigma \in G $ there is a natural auto-equivalence $ \tilde{\sigma}_{*} $ of the category
$ \DGMod(\AB) $
given by twisting the (right) action of $ \AB $ by $ \tilde{\sigma}^{-1} $. It is obvious that $ \tilde{\sigma}_{*} $
preserves the full subcategory of DG-modules satisfying \eqref{7.1}--\eqref{7.3}, and hence we have
a natural action of $G$ on the set $ \tilde{\MM} $
of isomorphism classes of such modules. On the other hand, $ \tilde{\sigma}_{*} $ preserves also
$ \Mod(A) $ (regarded as a subcategory of $ \DGMod(\AB) $) and, more specifically, the full
subcategory of rank $1$ torsion-free $A$-modules in $ \Mod(A)$. This gives an action of $ G $ on $ \RR$.
Now, if a quasi-isomorphism $\, q: M \to \LB \,$ satisfies Definition~\ref{D2} then obviously so does
$ \tilde{\sigma}_{*}q: \tilde{\sigma}_{*}M \to \tilde{\sigma}_{*}\LB \,$. Hence, the map
$\, \tilde{\omega}_1: \tilde{\MM} \to \RR \,$ defined by taking cohomology of DG-modules is
$G$-equivariant. By Theorem~\ref{T6} and Corollary~\ref{bij},  $\, \tilde{\omega}_1 \,$ is a
bijection equal to $ \omega_1 \circ \varrho_* \,$; the inverse map
$\, \tilde{\theta}_1 := \tilde{\omega}_1^{-1} =
\varrho_*^{-1} \circ \theta_1 \,$ is thus a $G$-equivariant bijection as well.

Next, by definition, $G$ is a subgroup of $ \Aut_k(R) \,$, so each $\,\sigma \in G\,$ gives a twisting
functor $ \sigma_{*} $ on the category $ \Com(R) $. The natural embedding
$\, \iota_{*}:\, \DGMod(\AB) \to \Com(R) \,$ intertwines $ \tilde{\sigma}_{*} $ and
$ \sigma_{*} $ for any $ \sigma \in G \,$. Hence, if we define the action of $ G $ on $ \CC $
by $\, [(\X, \Y, \i, \j)] \, \mapsto \, [(\sigma^{-1}(\X), \sigma^{-1}(\Y), \i, \j)] \, $
the map $\, \tilde{\theta}_2 := \theta_2 \circ \varrho_{*} : \tilde{\MM} \to \CC \,$
becomes $G$-equivariant. More precisely, by Theorem~\ref{T6} and Corollary~\ref{bij},
$ \tilde{\theta}_2 $
is a $G$-equivariant {\it bijection} and hence so is its inverse $\tilde{\omega}_2 $.

Thus, passing from $\A$-envelopes to DG-envelopes we can refine Theorem~\ref{T6}:
\begin{theorem}
\la{T13}
The Calogero-Moser correspondence factors through the $G$-equivariant bijective
maps:
\begin{equation}
\la{D8}
\begin{diagram}[small, tight]
\RR &  \pile{\rTo^{\tilde{\theta}_1}\\ \lTo_{\tilde{\omega}_1}} & \tilde{\MM} &
 \pile{\rTo^{\tilde{\theta}_2}\\ \lTo_{\tilde{\omega}_2}} & \CC \ , \\
\end{diagram}
\end{equation}
and hence is $G$-equivariant.
\end{theorem}
It remains to note that $ G $ is isomorphic to the automorphism group of the Weyl algebra (see \cite{M-L}),
and the actions of $G$ on $ \RR $ and $ \CC $ defined above are the same as in \cite{BW1, BW2}.

\section{Functoriality}

The classical (commutative) analogue of the Calogero-Moser correspondence relates the rank $1$ torsion-free
modules over the polynomial algebra $ k[x,y] $ to its finite-dimensional cyclic representations.
It is easy to see that this relation is functorial, the corresponding functor being the cokernel of the natural transformation $\,M \to M^{**}\,$. Such a construction, however, does not generalize immediately
to the noncommutative case since, unlike for $ k[x,y] $, all rank $1$ torsion-free modules over $ A_1 $ are projective and hence reflexive (meaning that $\, M \cong M^{**}\,$). To understand this apparent
``loss of functoriality'' was our original motivation for the present work.
We conclude the paper with a few remarks concerning this question. In fact, we  give
an answer (see Corollary~\ref{C6} and remarks thereafter), but probably not {\it the} answer, as some
more subtle questions seem to arise.

Passing to the category of $\A$-modules reduces the problem of functoriality to that of extension of
morphisms. Due to Theorem~\ref{T1} the last problem has a simple solution which can be stated as follows.
\begin{proposition}
\la{exte}
Let $\, r_1: M_1 \to \K_1 \,$ and $\, r_2: M_2 \to \K_2 \,$  be $\A$-envelopes of modules
$M_1$ and $ M_2$ respectively. Then every $A$-module homomorphism $\, f: M_1 \to M_2 \,$ extends
to a morphism $\, \tilde{f}: \K_1 \to \K_2 \,$ of $\A$-modules so that the diagram
\[
\begin{diagram}[small, tight]
M_1       &  \rTo^{r_1}    & \K_1        \\
\dTo^{f}  &                & \dTo_{\tilde{f}}  \\
M_2       & \rTo^{r_2}     & \K_2    \\
\end{diagram}
\]
commutes in $ \Modi(A)$. Such an extension is unique up to $\A$-homotopy.
\end{proposition}
\begin{proof}
By Theorem~\ref{T1}$(a)$, there exists an $\A$-morphism $\, s_1: \K_1 \to M_1 \,$ such that
$\, r_1 \circ s_1 \,$ is homotopy equivalent to the identity map $ \id_{\K} $ in $\, \Modi(A) \,$.
Given now a morphism $\, f: M_1 \to M_2 \,$ in $ \Mod(A) $ we set
$\,\tilde{f} := r_2 \circ f \circ s_1 \,$. The difference  $\, \tilde{f} \circ r_1 - r_2 \circ f \,$
is then a nullhomotopic map $ M_1 \to \K_2 $ in $ \Modi(A) $, and hence is zero by degree considerations.
If $ \tilde{f}': \K_1 \to \K_2 $ is another extension of $ f $ we have $\,
(\tilde{f}' - \tilde{f}) \circ r_1 = 0 $, and hence $\, (\tilde{f}' - \tilde{f})\circ r_1 \circ s_1 = 0 \, \Rightarrow \, \tilde{f}' - \tilde{f} \sim 0 \,$ in $ \Modi(A)$.
\end{proof}
\begin{corollary}
\la{C6}
Choosing an $\A$-envelope, one for each $\, M \,$, and assigning to each homomorphism
$\, f: M_1 \to M_2  \,$ the $\A$-homotopy class $\, [\,\tilde{f}\,]_{\infty} \,$ of its extension $ \tilde{f} $
defines an equivalence $\, \Theta \, $ from the full subcategory $\, \Ideals(A) \,$ of
rank $ 1 $ torsion-free modules in $ \Mod(A) $
to the full subcategory of $ \Di(A) $ consisting of $ \A$-modules satisfying
\eqref{2.1}--\eqref{2.5}.
\end{corollary}
\begin{proof} By Proposition~\ref{exte}, $\, \Theta \,$ is a well-defined functor.
Its inverse is given by taking cohomology of an $\A$-envelope.
\end{proof}

\noindent
We close this section with a few general remarks.

1. Combining Corollary~\ref{C6} with Theorem~\ref{T3} we see that the bijection
$\, \RR \to \MM \,$ of Theorem~\ref{bij} is induced by an equivalence $\Theta $.
This recovers, at least in part, the functoriality of the Calogero-Moser correspondence
in the noncommutative case.

2. If $\, A = k[x,y] \,$, every $\A$-module satisfying \eqref{2.1}--\eqref{2.5} is a
strict complex of $A$-modules (see Proposition~\ref{P1}). Moreover, in this case every
$A$-module map $\, f: M_1 \to M_2 \,$ extends to a unique $A$-linear morphism
$\, \tilde{f}: \K_1 \to \K_2 \,$, and the equivalence $\, \Theta \, $ of Corollary~\ref{C6}
factors thus through $ \Com(A)$:
$$
\Ideals(A) \stackrel{\Theta_0}{\longrightarrow} \Com(A) \stackrel{\Upsilon}{\longrightarrow}
\Modi(A) \to \Di(A)\ ,
$$
$\, \Theta_0 \,$ being the cokernel functor
$\, \Theta_0(M) := (\,M^{**} \onto M^{**}/M \,)\,$ mentioned above.

3. The preceeding remark shows that in the commutative case every homotopy class of extensions
$ [\,\tilde{f}\,]_{\infty} $ contains a unique {\it strict} representative. One might wonder if this is true for the Weyl algebra. By Theorem~\ref{T3}, every {\it iso}morphism indeed extends to a strict isomorphism of
$\A$-envelopes but, in general, the answer is negative. In fact, if every $\, f \in \End_{A}(M) \,$ were
extendable to a strict endomorphism $\, \tilde{f}: \K \to \K \,$, we would have a
non-trivial representation of $ \End_{A}(M) $ on the finite-dimensional vector space $ K^1 $.
But $ \End_{A}(M) $ is Morita equivalent to $ A $ and hence cannot have nonzero finite-dimensional modules.
In view of Theorem~\ref{T12}, this implies also that Proposition~\ref{exte} does not hold for DG-envelopes:
to be precise, {\it not} every $A$-module map $\, f: M_1 \to M_2 \,$ extends through a DG-envelope $\,
M_1 \to \L_1 \,$ to a morphism $\, \L_1 \to \L_2 \,$ in $ \DGMod(\AB) $.

4. It is still an interesting question whether there exist some
distinguished non-strict extensions of morphisms $\,f \,$ in $\Mod(A)$. 
We expect that there is a natural functor (in fact, an $\A$-functor) 
$\, \tilde{\Theta}: \,\Ideals(A) \to \MOD_{\infty}(A) \,$, which takes 
values in a {\it DG-category} of $\A$-modules (see \cite{D}, \cite{K1}, \cite{Ko})
and descends to $\, \Theta \,$ at the cohomology level.

\section{Appendix: A DG-structure on Local Cohomology}
\la{App}
In this appendix we will give a geometric construction of DG-envelopes using the language of noncommutative
projective schemes (see \cite{AZ}, \cite{SV}). This construction is less elementary than the ones described
in Section~\ref{Ex}; however, apart from clarifying the relation to the previous work \cite{BW2}, it gives
a cohomological interpretation of DG-envelopes and exhibits a geometric origin of the properties axiomatized
in Definition~\ref{D2}.

\subsection{Projective closure}
In algebraic geometry, there is a standard procedure of passing from affine varieties to projective ones:
given an affine variety with a fixed embedding in an affine space, say $\, X \subseteq \AA^{n}_k \,$, one identifies $ \AA^{n}_k $ with the open complement of a coordinate hyperplane in $ \PP^n_k $ and takes
the closure of $ X $ in $ \PP^n_k $. In this way one gets a projective variety $ \bX $
containing $ X $ as an open subset whose complement $\, Z = \bX \setminus X \,$ is an ample divisor
in $ \bX $.

This procedure generalizes to the realm of noncommutative geometry as follows
(see, e.g., \cite{LVdB}, \cite{LeB}, \cite{Sm}).
Let $A$ be a finitely generated associative $k$-algebra. If we think of the category of
noncommutative affine schemes over $k$ as the dual to that of associative $k$-algebras, the free algebra
$\, R = k \langle x_1, x_2, \ldots, x_n \rangle \,$ corresponds to the noncommutative affine space
$\,  N\AA^{n}_{k} \,$ (cf. \cite{KR}), while an epimorphism of algebras $\,R \onto A\,$ to a closed
embedding $ X \into N\AA^{n}_{k} $. The natural filtration on $ R $ (defined by giving each generator
$ x_i $ degree $1$) descends to a positive filtration
$ \{A_{\bullet}\} $  on $ A $, and we may form the graded Rees algebra
$$
\tA := \bigoplus_{i \geq 0} A_i  t^i \subset A[t] \ .
$$
The projective closure of $ X $ is then defined categorically, in terms of graded
$ \tA$-modules. More precisely, we identify $ \bX $ with the category $ \Qcoh(\bX) $ of
quasicoherent sheaves on it, which, in turn, is defined as the quotient category
$\, \Qgr(\tA)  \,$ of (right) graded $\tA$-modules modulo torsion (see \cite{AZ}).
Thus, by definition, $\,\Qcoh(\bX) \equiv \Qgr(\tA) \,$ is a $k$-linear Abelian category
which has enough injectives and comes equipped with two natural functors --- the exact
quotient functor $\, \pi:  \GrM(\tA) \to \Qcoh(\bX)\,$ and its right adjoint (and hence,
left exact) functor $\, \omega: \, \Qcoh(\bX) \to  \GrM(\tA) \,$. 
We introduce the {\it cohomology of quasicoherent sheaves on $\bX$} by taking the right derived 
functors of $\,\omega\,$:
\begin{equation}
\la{AA1}
\bH^n(\bX,\,\ms{F}) := {\mathsf R}^{n}\omega (\ms{F}) \ ,  \quad \ms{F} \in \Qcoh(\bX) \ .
\end{equation}
Then $\, \bH^n(\bX, \ms{F}) $ has a natural structure of graded $\tA$-module: its $m$-th graded
component can be identified with $\, {H}^n(\bX, \ms{F}(m)) := \Ext^{n}(\pi\tA,\, \ms{F}(m))\,$,
where $\, \ms{F}(m) \in \Qcoh(\bX) \,$ is a ``twisted'' sheaf obtained from $ \ms{F} $ by
shifting the grading of $ \tA$-modules by the integer $\, m \in \Z \,$.

Setting $\,\Qcoh(X) := \Mod(A) \,$ and $\,\Qcoh(Z) := \Qgr(\tA/(t)) \,$, we think 
of $ X $ geometrically as an open affine subvariety of $ \bX $ and of $Z$ as the divisor at 
infinity that complements $ X $ in $ \bX $. In this situation we have the diagram of standard 
functors:
$$
\begin{diagram}[small, tight]
 \Qcoh(X) & \ \pile{\lTo^{\,\,j^*}\\ \rTo^{\,\,j_*}} \ &  \Qcoh(\bX)  & \ \pile{\rTo^{i^*}\\
\lTo^{i_*}\\ \rTo^{i^!}} \ & \Qcoh(Z)   \\
\end{diagram}
$$
imitating the usual relation between sheaves on $ \bX $ and its open and closed
subspaces  $\, X \stackrel{j}{\into} \bX \stackrel{i}{\otni} Z\,$ (see \cite{Sm}, Section~8).
\begin{lemma}
\la{LA1}
For every $\, \ms{F} \in \Qcoh(\bX) \,$ we have the exact sequence in $ \Qcoh(\bX) $
\begin{equation}
\la{A2}
0 \to i_* i^{!} \ms{F}(-1) \to \ms{F}(-1) \stackrel{\cdot\,t}{\longrightarrow} \ms{F} \to i_* i^* \ms{F} \to 0\
\end{equation}
with map in the middle induced by multiplication by $ t \in \tA$.
\end{lemma}
\begin{proof}
At the level of graded modules, the functors $\, i^*,\, i^{!}:\, \GrM(\tA) \to \GrM(\tA/(t)) \,$
are defined by
$$
i^* \tM := \tM \otimes_{\tA} \tA/(t) \cong \tM/\tM t \quad ,
\quad i^! \tM := \Ker\,[\tM \stackrel{\cdot t}{\longrightarrow} \tM(1)] \ ,
$$
so for every $ \tM $ we have the exact sequence in $\,\GrM(\tA) \,$:
\begin{equation}
\la{AA2}
0 \to i_* i^{!} \tM(-1) \to \tM(-1) \stackrel{\cdot\,t}{\longrightarrow} \tM \to i_* i^* \tM \to 0\ .
\end{equation}
Since each $ \ms{F} \in \Qcoh(\bX) $ can be represented by a graded module, $\, \ms{F} = \pi\,\tM \,$,
and the quotient functor $\,\pi: \, \GrM(\tA) \to \Qgr(\tA/(t)) \,$ is exact,
\eqref{AA2} implies  \eqref{A2}.
\end{proof}
\subsection{Cohomology with supports}
Given a graded module $\, \tM \in \GrM(\tA) \,$, we write $\, \ga \tM \in \GrM(\tA) \,$ for its largest
$t$-torsion submodule:
$$
\ga \tM = \{ m \in \tM \,: \, m\,t^n = 0 \ \mbox{for some}\ n\geq 0 \}\ .
$$
If $\, f: \tM \to \tN \, $ is a morphism in
$\,\GrM(\tA)\,$, we have $\, f(\ga\tM) \subseteq \ga \tN \,$, so there is a map $ \ga(f):
\ga\tM \to \ga \tN $ which agrees with $ f $ on each element of $ \ga \tM\,$. Thus
$\, \ga: \, \GrM\,\tA \to \GrM\,\tA $ is an additive functor on graded modules. Equivalently,
we can define $ \ga $ by the formula
\begin{equation}
\la{A3}
\ga \tM = \lim\limits_{\longrightarrow}^{}\,\bHom(\tA/(t)^n,\, \tM) \ ,
\end{equation}
where $\, \bHom \,$ stands for the graded functor describing $\tA$-module homomorphisms
of all (finite) degrees, and the system of $ \bHom$'s is directed by restrictions of such
homomorphisms through the canonical algebra maps $\, \tA/(t)^{n+1} \onto \tA/(t)^n$. Since
$ \, \lim\limits_{\longrightarrow}^{} $ is exact and each $\, \bHom(\tA/(t)^n,\,\mbox{---}\,) \,$
is left exact, it follows from \eqref{A3} that $ \ga $ is a left exact functor.

Now, given $\, \ms{F} \in \Qcoh(\bX) \,$, we set $\, \bH^{0}_{\!Z}(\bX,\,\ms{F}) := \ga \,\omega \ms{F}\,$
and think of $\,\bH^{0}_{\!Z}(\bX,\,\ms{F}) \, $ geometrically as the space of sections of the twisted sheaf
$\, \bigoplus_{m \in \Z}\,\ms{F}(m) \,$ supported on the divisor $Z$.
By definition, $\, \bH^{0}_{\!Z} \,$ is the composite of two left exact functors, and hence
left exact. We define the higher {\it cohomology of sheaves with support in Z} as the right derived 
functors of $ \bH^{0}_{\!Z}\,$,
i.~e.
\begin{equation}
\la{A4}
\bH_{\!Z}^n(\bX,\,\ms{F}) := {\mathsf R}^{n}(\ga \,\omega) \,\ms{F} \ , \quad n \geq 0 \ .
\end{equation}
For each $ n \geq 0\,$, $\,\bH_{\!Z}^{n}(\bX,\,\ms{F})\,$ is then a graded $ \tA$-module, and we write
$$
\bH_{\!Z}^n(\bX, \ms{F}) := \bigoplus_{m \in \Z} {H}_{\!Z}^n(\bX, \ms{F}(m)) \ ,
$$
with $\, {H}_{\!Z}^{n}(\bX,\,\ms{F}(m))\,$ denoting the $m$-th graded component of
$\,\bH_{\!Z}^{n}(\bX,\,\ms{F})\,$.

\begin{proposition}
\la{PA1}
For every $\, \ms{F} \in \Qcoh(\bX) \,$  there is an exact sequence in $\,\GrM(\tA)$
$$
0 \to \bH^{0}_{\!Z}(\bX,\,\ms{F}) \to \bH^{0}(\bX,\,\ms{F}) \to \bH^{0}(X,\, j^*\ms{F})
\stackrel{q}{\longrightarrow} \bH^{1}_{\!Z}(\bX,\,\ms{F}) \to \bH^{1}(\bX,\,\ms{F}) \to 0
$$
and isomorphisms $\,\bH^{n}_{\!Z}(\bX,\,\ms{F}) \cong \bH^{n}(\bX,\,\ms{F})\,$ for all
$\, n \geq 2\,$.
\end{proposition}
\begin{proof}
The standard proof of this result in the geometric case (see \cite{H}, Corollary~1.9)
involves flasque resolutions of sheaves which are not defined in our categorical setting.
We will use instead injective resolutions which are available, since the quotient category
$ \Qcoh(\bX)\,$ has enough injectives (see, e.g., \cite{AZ}, Proposition~7.1).

Thus, let $\, \ms{F} \to \ms{E}^{\bullet} \,$ be an injective resolution of $ \ms{F} $ in $ \Qcoh(\bX)\,$.
For each $ n \geq 0 \,$, set $\, \tE^n := \omega \ms{E}^{n} \,$, $\, \tI^{n} := \ga \tE^{n}
\subseteq \tE^n \,$ and $\, \tQ^n :=  \tE^n/\tI^n\,$. Then there is an exact sequences of
complexes in $ \, \GrM(\tA)\,$:
\begin{equation}
\la{A5}
0 \to \tI^{\bullet} \to \tE^{\bullet} \to \tQ^{\bullet} \to 0\ ,
\end{equation}
which gives a long exact sequence in cohomology:
\begin{equation}
\la{A6}
\ldots \,\to H^{n}(\tI^{\bullet}) \to H^{n}(\tE^{\bullet}) \to
H^{n}(\tQ^{\bullet}) \to H^{n+1}(\tI^{\bullet}) \to \, \ldots
\end{equation}
With definitions \eqref{AA1} and \eqref{A4}, we have at once
$\, H^{n}(\tE^{\bullet}) \cong  \bH^{n}(\bX,\,\ms{F})\,$  and
$\, H^{n}(\tI^{\bullet}) \cong  \bH_{\!Z}^{n}(\bX,\,\ms{F}) \,$
for all $ n \geq 0 \,$. On the other hand, the functors $ j^* $ and $ j_* $
are both exact and send injectives to injectives. Hence $\, j_*j^*\ms{F} \to
j_*j^*\ms{E}^{\bullet} \,$ is an injective resolution of $ j_*j^*\ms{F} $,
and we have $\, \bH^{n}(\bX,\, j_* j^*\ms{F})  \cong ({\mathsf R}^n \omega)\,j_* j^*\ms{F}
\cong H^n(\omega j_* j^*\ms{E}^{\bullet})\,$.

Now, by definition, $ j^* $ and $ j_* $ factor through $\,\GrM(\tA) \,$,
i.~e. $\, j^* \cong \tilde{j}^* \omega \,$ and $\, j_* \cong \pi \tilde{j}_*\,$,
where $\, \tilde{j}^*:\, \GrM(\tA) \to \GrM(\tA[t^{-1}])\,$ is the (graded) localization functor
and $\, \tilde{j}_* \,$ its right adjoint (the restriction functor).
Hence we have $\, \omega j_* j^*\ms{E}^\bullet \cong \omega \pi \tilde{j}_* \tilde{j}^* \omega \ms{E}^\bullet \cong \omega \pi \tilde{j}_* \tilde{j}^* \tE^\bullet $. Since $\, \tilde{j}^* \,$ is exact and
$\, \tilde{j}^* \tI^n = 0 \,$ for all $\,n \geq 0\,$, it follows from \eqref{A5} that
$ \tilde{j}^* \tE^{\bullet} \cong \tilde{j}^* \tQ^{\bullet} $. By construction, each $ \tQ^n $ is a
$t$-torsion-free (and hence, torsion-free) injective module, so
$\, \tilde{j}_* \tilde{j}^* \tQ^{n} \cong \tQ^n \,$ and $\, \omega \pi \tQ^n \cong \tQ^n \,$.
Combining the above isomorphisms together, we get
$$
\omega j_* j^*\ms{E}^\bullet \,\cong\,
\omega \pi \tilde{j}_* \tilde{j}^* \tE^{\bullet} \,\cong\, \omega \pi
\tilde{j}_* \tilde{j}^* \tQ^{\bullet}\, \cong \,\omega \pi \tQ^\bullet \,\cong \,\tQ^\bullet \ ,
$$
and thus $\,\bH^{n}(\bX,\, j_* j^*\ms{F})  \cong H^n(\tQ^\bullet) \,$ for all $ n \geq 0 $.
The long cohomology sequence \eqref{A6} now becomes
\begin{equation}
\la{AA6}
\ldots \, \to \bH^{n}_{\!Z}(\bX,\,\ms{F}) \to \bH^{n}(\bX,\,\ms{F}) \to \bH^{n}(\bX,\, j_*j^*\ms{F})
\to \bH^{n+1}_{\!Z}(\bX,\,\ms{F}) \to \, \ldots
\end{equation}
To finish the proof it remains to show that $\, \bH^{0}(\bX,\, j_*j^*\ms{F}) \cong
\bH^{0}(X,\, j^*\ms{F}) \,$ and $\, \bH^{n}(\bX,\, j_*j^*\ms{F}) = 0 \,$ for
$ n\geq 1 $. This follows at once from the definition of $\, \bH^{n}(X,\, j^*\ms{F}) \,$:
\begin{equation}
\la{A7}
\bH^{n}(X,\, j^*\ms{F}) := \bigoplus_{m \in \Z} {\rm Ext}_{A}^{n}(A,\,j^*\ms{F}(m)) \cong
\left\{
\begin{array}{cl}
\oplus_{m \in \Z}\ j^*\ms{F}  & \mbox{if}\ n = 0 \\*[1ex]
0 &  \mbox{otherwise}
\end{array}
\right.
\end{equation}
and the Leray type spectral sequence (cf. \cite{Sm}, Lemma~2.8):
\begin{equation}
\la{A8}
{\bH}^p(\bX,\, {\mathsf R}^qj_*j^* \ms{F}) \ \Rightarrow \ {\bH}^{p+q}(X,\, j^*\ms{F})\ .
\end{equation}
In view of exactness of $ j_* $, the sequence \eqref{A8} collapses on the line $ q = 0 $
giving the required isomorphisms  $\, \bH^{n}(\bX,\, j_*j^*\ms{F}) \cong \bH^{n}(X,\, j^*\ms{F}) \,$
for all $ n \geq 0 $.
\end{proof}
\subsection{Envelopes} Now we are in position to formulate conditions underlying
the existence of DG-envelopes in the geometric setting.

First, recall that we are working with an algebra $ A $ with a fixed set of generators
$ \{x_1, x_2, \ldots, x_n\} $, or equivalently, with an algebra epimorphism
$\, \eta:\, R \onto A \,$, where $\, R = k \langle x_1, x_2, \ldots, x_n\rangle$.
As in Section~\ref{HPT}, we form the two-component DG-algebra $\, \AB = I \oplus R \,$ with
$\, I = \Ker(\eta) \,$ and differential given by the natural inclusion $ d:\, I \into R $ .
The map $ \eta $ extends then to a quasi-isomorphism of DG-algebras $\, \AB \to A \,$ which we
will also denote by $\eta$.

Second, we need some finiteness results, and thus, will work with objects
$ \ms{F} \in \Qcoh(\bX) $ represented by finitely generated $\tA$-modules.
These form a full subcategory of $ \Qcoh(\bX) $ which, by analogy with the geometric case, is called
the category $\,\Coh(\bX) \,$ of {\it coherent sheaves} on $ \bX $. In addition, we assume that
$ \tA $ is Noetherian and satisfies the Artin-Zhang property $ \chi $ (see \cite{AZ}, Definition~3.7):
this guarantees that Serre's Finiteness and Vanishing theorems
hold for the cohomology of coherent sheaves on $ \bX $ (see \cite{AZ}, Theorem~7.4).

Now, let $\, M \in \Mod(A) \,$ be a finitely generated $A$-module. We call $\, \ms{F} \in \Coh(\bX)\,$ an {\it extension} of $ M $ to $ \bX $ if $\, j^*\ms{F} = M \,$ and $\, i^{!}\ms{F} = 0\,$\footnote{The last
condition means that $ \ms{F} $ arises from the ``homogenization'' of $ M $ relative to some module
filtration.}. The following theorem is the main result of this section.
\begin{theorem}
\la{TA1}
Assume that $\, \ms{F} \in \Coh(\bX) \,$ is an extension of $ M $ to $ \bX $ satisfying
 $\,(a)\ H^0(\bX,\, \ms{F}) = 0 \,$ and $\,(b)\ H^0(Z,\, i^*\ms{F}) = H^{1}(Z,\, i^*\ms{F}) = 0\,$.
Then the complex of vector spaces
\begin{equation}
\la{A9}
\L := [\,0 \to H_{\!Z}^1(\bX,\, \ms{F}) \to H^1(\bX,\, \ms{F}) \to 0\,]
\end{equation}
has a natural structure of DG-module over $ \AB $ with connecting morphism $\, q\,$
(see Proposition~\ref{PA1}) giving a quasi-isomorphism $\, M \to \L \,$ in $\, \DGMod(\AB) $.
\end{theorem}
\begin{proof}
We combine Lemma~\ref{LA1} and Proposition~\ref{PA1} to get the following commutative diagram
with exact rows and columns:
\begin{equation}
\la{DA2}
\begin{diagram}[small, tight]
  &      &  0   & \rTo & H_{\!Z}^0(\bX,\, i_*i^*\ms{F}) & \rTo & H^0(\bX,\, i_*i^*\ms{F}) & \rTo & 0 \\
  &      & \dTo &      &           \dTo                 &      &       \dTo               &      &   \\
0 & \rTo &  M   & \rTo &  H_{\!Z}^1(\bX,\, \ms{F}(-1))  & \rTo &  H^1(\bX,\, \ms{F}(-1))  & \rTo & 0 \\
  &      & \dEq &      &        \dTo_{\cdot\,t}         &      &  \dTo_{\cdot \, t}       &      &   \\
0 & \rTo &  M   & \rTo &   H_{\!Z}^1(\bX,\, \ms{F})     & \rTo &   H^1(\bX,\, \ms{F})     & \rTo & 0 \\
  &      & \dTo &      &           \dTo                 &      &       \dTo               &      &   \\
  &      &  0   & \rTo & H_{\!Z}^1(\bX,\, i_*i^*\ms{F}) & \rTo & H^1(\bX,\, i_*i^*\ms{F}) & \rTo & 0 \\
\end{diagram}
\end{equation}
To be precise, the two rows in the middle are the degree $-1$ and $0$ components of
the 5-term exact sequence of Proposition~\ref{PA1}. Each of these begins with $\,0\,$ because
of the vanishing of $ H^0(\bX,\, \ms{F}) $ and $ H^0(\bX,\, \ms{F}(-1)) $.
(Note that $\, H^0(\bX,\, \ms{F}) = 0 \,  \Rightarrow \, H^0(\bX,\, \ms{F}(-1)) = 0 \,$
by Lemma~\ref{LA1}, since $\, i^! \ms{F} = 0 \,$.) The spaces $ H^0(X, \, j^*\ms{F}(-1)) $ and
$ H^0(X, \, j^*\ms{F}) $ can then be identified with $ j^*\ms{F} = M $ by \eqref{A7}.
The first and the last rows are also part of the exact sequence of Proposition~\ref{PA1}
with $ \ms{F} $ being replaced by $ \, i_*i^* \ms{F} \,$.
In this case, the sequence in question breaks up into two isomorphisms due to the natural
identity $\, j^* i_* = 0 \,$ (see \cite{Sm}, Proposition~8.5). Finally, the columns of \eqref{DA2}
arise from applying the functors $ H^0_{\!Z} $ and $ H^0 $ to the short exact sequence
\eqref{A2} (with first term vanishing), and thus are also exact.

Now, according to \cite{AZ}, Theorem~8.3, $\, H^n(\bX,\, i_*i^*\ms{F}) \cong
H^n(Z,\, i^*\ms{F}) \,$ for all $ n \geq 0 $. Hence, with our assumptions on $ \ms{F} $,
the first and the last rows of \eqref{DA2} vanish, and the map induced by multiplication
by $t$ is an isomorphism. If we identify the complex
$\, 0 \to H_{\!Z}^1(\bX,\, \ms{F}(-1)) \to H^1(\bX,\, \ms{F}(-1)) \to 0 \,$
with $ \L $ via this isomorphism, then $ \L $ gets naturally a structure of $R$-module. Indeed,
each of the generators $\, \tx_1 = x_1 t,\,\ldots, \, \tx_n = x_n t\, $ of $ \tA $
has degree $ 1 $ and hence, under our identification, induces a linear endomorphism of $ \L $.
As $ R $ is free, we get thus a homomorphism of algebras $\,
\alpha:\, R \to \End_{k}(\L)^{\rm opp} $ defining a right action of $ R $ on $ \L $.
When restricted to cohomology, this action coincides with the given action of $ A $ on $ M $.
Therefore, letting $\,(u, \,v).a := (- d_{L}^{-1}(v) \,d_{A}(a), \,0)\,$ for
$ a \in I $ and $ (u,\,v) \in \L \,$, we may extend  $\, \alpha \,$ to a (unique) map of DG-algebras:
$ \AB \to \underline{\End}_k^\bullet(\L)^{\rm opp} $. This gives $ \L $ a structure of
DG-module over $ \AB $, with connecting morphism $\, M \cong H^0(X,\, j^*\ms{F})
\stackrel{q}{\longrightarrow} \L \,$ becoming a quasi-isomorphism in $ \DGMod(\AB)\,$.
\end{proof}
\subsection{The Weyl algebra} Now we return to our basic example: thus,
let $ A $ be again the Weyl algebra with canonical generators $x$ and $y$,
$ \, \eta: R \onto A $, the corresponding projection from
$\, R = k\langle x,y\rangle \,$, and $\, \AB = RwR \oplus R \,$, the DG-algebra defined in Section~\ref{HPT}.
In this case, $\tA$ is generated by the elements
$\, \tx = x t \,$, $\, \ty = y t \,$ and $\, t \,$, all having degree $1$
and satisfying the relations $\, [\tx, \,t] = [\ty,\,t] = 0 \,$ and $\, [\tx,\,\ty] = t^2 $.
The closure of $ X $ is then a quantum projective plane $ \PP^2_{\!q} $ (see \cite{A}), and
$ Z = \bX \setminus X $ is the usual (commutative) projective line $ \PP^1 $.
Now, according to \cite{BW2}, Lemma~4.1, every f.~g. rank $1$ torsion-free $A$-module $ M $ has a (unique)
extension $\, \ms{M} \,$ to $\, \PP^2_{\!q} \,$ such that $\, i^* \ms{M} \cong \ms{O}_{\PP^1} \,$.
To apply Theorem~\ref{TA1}, we set $\, \ms{F} := \ms{M}(-1) \,$. The condition $(a)$ follows
then from \cite{BW2}, Theorem~4.5(ii), while $ (b) $ holds automatically, since $ \ms{O}_{\PP^1}(-1) $ is an
acyclic line bundle on $ \PP^1 $. Thus, the complex $ \L $ given by \eqref{A9} is a DG-module over
$ \AB $, and there is a quasi-isomorphism $\, q:\, M \to \L \,$ in $ \DGMod(\AB)$.

We claim that $\, \L \,$ is a DG-envelope of $ M $ in the sense of Definition~\ref{D2}. Indeed,
the axiom \eqref{7.1} follows at once from Serre's Finiteness theorem.
To verify \eqref{7.2} and \eqref{7.3}, observe first that $\, \bH^{0}_{\!Z}(\PP^2_{\!q},\,\ms{M}) = 0\,$
as $ \bH^{0}(\PP^2_{\!q},\,\ms{M}) $ is torsion-free (and hence, $t$-torsion-free) by \cite{BW2}, Proposition~4.3.
Applying now $\, H^{0}_{\!Z} \,$ to the short exact sequence
$$
0 \to \ms{M}(-1) \to \ms{M} \to i_* \ms{O}_{\PP^1} \to 0
$$
and taking into account the isomorphisms
$$
H^{n}_{\!Z}(\PP^2_{\!q},\, i_* \ms{O}_{\PP^1}) \cong H^{n}(\PP^2_{\!q}, \, i_*\ms{O}_{\PP^1})
\cong H^{n}(\PP^1, \, \ms{O}_{\PP^1})  \quad \mbox{for}\ n = 0, 1 \quad \mbox{(see \eqref{DA2})}
$$
we get the exact sequence
$$
0 \to H^{0}(\PP^1,\, \ms{O}_{\PP^1}) \stackrel{\delta}{\longrightarrow}
H^{1}_{\!Z}(\PP^2_{\!q},\, \ms{M}(-1)) \stackrel{t}{\longrightarrow} H^{1}_{\!Z}(\PP^2_{\!q},\, \ms{M}) \to
H^{1}(\PP^1, \, \ms{O}_{\PP^1})\ .
$$
By Liouville's Theorem, we have $\, H^{0}(\PP^1, \, \ms{O}_{\PP^1}) \cong k \,$, while
$ H^{1}(\PP^1, \, \ms{O}_{\PP^1}) = 0 \,$. Now, let
$\, i \in L^0 := H_{\!Z}^1(\PP^2_{\!q},\,\ms{M}(-1))\,$ be the image of a basis vector of
$\, H^{0}(\PP^1, \, \ms{O}_{\PP^1}) \,$ under the connecting map $ \delta $. Then $\, i \,$ spans
the kernel of $\, t\,$, and the fact that it is a cyclic generator of the $R$-module $ L^0 $ follows
from \cite{BW2}, Lemma~6.1. Writing $\,X \,$ and $\, Y \,$ for
the action of $ \tx $ and $ \ty $  on $\, H_{\!Z}^1(\PP^2_{\!q},\,\ms{M}(-1)) \,$,
we have $\, t \, (XY - YX + \id) = 0 \,$.
On the other hand, by our definition of the DG-structure on $ \L $,
$\, v.w = (XY - YX + \id)\,d^{-1}_L(v)\,$ for all $\, v  \in L^1 := H^1(\PP^2_{\!q},\,\ms{M}(-1)) $.
Whence $\, L^1.w \subseteq \Ker(t) = k.i \,$, which is equivalent to the axiom \eqref{7.3}.

\end{document}